\documentclass[a4paper,12pt]{amsart}
\usepackage{amsmath}
\usepackage{amssymb}
\usepackage{amsthm}
\usepackage[all]{xy}
\usepackage[latin1]{inputenc}        
\usepackage[dvips]{graphics}
\usepackage[dvips]{graphicx}
\usepackage{amscd}
\usepackage{fullpage}
\author{Francesco Polizzi}
\address{Dipartimento di Matematica, Universit{\`a} di Roma ``Tor Vergata'', Via della Ricerca Scientifica 1, 00133 Roma (Italy)}
\email{polizzi@.mat.uniroma2.it}

\title[Surfaces with $p_g=q=1, \; K^2=8$]{Surfaces of general type with $p_g=q=1, \; K^2=8$ and bicanonical map of degree $2$}
\date{\today}

\newtheorem{inizio}{Lemma}[section]
\newtheorem{theorem}[inizio]{Theorem}
\newtheorem{corollary}[inizio]{Corollary}
\newtheorem{proposition}[inizio]{Proposition}
\newtheorem{lemma}[inizio]{Lemma} 
\newtheorem{claim}[inizio]{Claim} 
\newtheorem{definition}[inizio]{Definition}
\newtheorem{remark}[inizio]{Remark}

\newcommand{\hB}{\widehat{B}}
\newcommand{\hS}{\widehat{S}}
\newcommand{\hmL}{\widehat{\mathcal{L}}}
\newcommand{\mO}{\mathcal{O}}
\newcommand{\hR}{\widehat{R}}
\newcommand{\hW}{\widehat{W}}
\begin{document}
\abstract
We classify the minimal algebraic surfaces of general type with $p_g=q=1, \; K^2=8$ and bicanonical map of degree $2$. It will turn out that they are isogenous to a product of curves, i.e. if $S$ is such a surface, then there exist two smooth curves $C, \; F$ and a finite group $G$ acting freely on $C \times F$ such that $S = (C \times F)/G$. We describe the $C, \; F$ and $G$ that occur. In particular the curve $C$ is a hyperelliptic-bielliptic curve of genus $3$, and the bicanonical map $\phi$ of $S$ is composed with the involution $\sigma$ induced on $S$ by $\tau \times id: C \times F  \longrightarrow C \times F$, where $\tau$ is the hyperelliptic involution of $C$. In this way we obtain three families of surfaces with $p_g=q=1, \; K^2=8$ which yield the first-known examples of surfaces with these invariants. We compute their dimension and we show that they are three generically smooth, irreducible components of the moduli space $\mathcal{M}$ of surfaces with $p_g=q=1, \; K^2=8$. Moreover, we give an alternative description of these surfaces as double covers of the plane, recovering a construction proposed by Du Val.    
\endabstract

\maketitle


\section{Introduction}
In [Par03] R. Pardini classified the minimal surfaces $S$ of general type with $p_g=q=0, \;  K_S^2=8$ and a rational involution, i.e. an involution $\sigma: S \longrightarrow S$ such that the quotient $T:=S / \sigma$ is a rational surface. All the examples constructed by Pardini are $\emph{isogenous to a product}$, i.e. there exist two smooth curves $C, \; F$ and a finite group $G$ acting faithfully on $C, \; F$ and whose diagonal action is free on the product $C \times F$, in such a way that $S = (C \times F)/G$. Pardini's classification contains five families of such surfaces; in particular, four of them are irreducible components of the moduli space of surfaces with $p_g=q=0, \; K_S^2=8$, and represent the surfaces with the above invariants and non-birational bicanonical map. \\
In this paper we deal with the irregular case, in fact we study the case $p_g=q=1, \; K_S^2=8$. Surfaces with $p_g=q=1$ are the minimal irregular surfaces of general type with the lowest geometric genus, therefore it would be very interesting to obtain their complete classification; for such a reason, they are currently an active topic of research. However, such surfaces are still quite mysterious, and only a few families have been hitherto discovered. If $S$ is a surface with $p_g=q=1$, then $2 \leq K_S^2 \leq 9$; the case $K_S^2=2$ is studied in [Ca81], whereas [CaCi91] and [CaCi93] deal with the case $K_S^2=3$. For higher values of $K_S^2$ only some sporadic examples were so far known; see [Ca99], where a surface with $K_S^2=4$ and one with $K_S^2=5$ are constructed. \\
When $p_g=q=1$, there are two basic tools that one can use in order to study the geometry of $S$: the \emph{Albanese fibration} and the \emph{paracanonical system}.  First of all, $q=1$ implies that the Albanese variety of $S$ is an elliptic curve $E$, hence the Albanese map $\alpha: S \longrightarrow E$ is a connected fibration; we will denote by $F$ the general fibre of $\alpha$ and by $g=g(F)$ its genus. Let us fix a zero point $0 \in E$, and for any $t \in E$ let us write $K_S+t$ for the line bundle $K_S+F_t-F_0$.  By Riemann-Roch and semicontinuity theorem we have $h^0(S, \;K_S +t)=1$ for general $t \in E$, hence denoting by $C_t$ the only element in the complete linear system $|K_S+t|$ we obtain a $1-$dimensional algebraic family $\{K \}= \{C_t\}_{t \in E}$ parametrized by the elliptic curve $E$. We will call it the paracanonical system of $S$; according to [Be88], it is the irreducible component of the Hilbert scheme of curves on $S$ algebraically equivalent to $K_S$ which dominates $E$. The \emph{index} $\iota=\iota(K)$ of the paracanonical system $\{K \}$ is the number of distinct curves of $\{K \}$ through a general point of $S$. The \emph{paracanonical map} $\omega: S \longrightarrow E(\iota),$ where $E(\iota):=\textrm{Sym}^{\iota}E$, is defined in the following way: if $x \in S$ is a general point, then $\omega(x)=t_1+ \cdots +t_{\iota}$, where $C_{t_1}, \ldots, C_{t_{\iota}}$ are the paracanonical curves containing $x$. The best result that one might obtain would be to classify the triples $(K^2,g, \iota)$ such that there exists a minimal surface of general type $S$ with $p_g=q=1$ and these invariants. Since by the results of Gieseker the moduli space $\mathcal{M}_{\chi, \; K^2}$ of surfaces of general type with fixed $\chi(\mathcal{O}_S), \; K_S^2$ is a quasiprojective variety, it turns out that there exist only finitely many such triples, but a complete classification is still missing. \\By the results of [Re88], [Fr91] and [CaCi91] it follows that the bicanonical system $|2K_S|$ of a minimal surface of general type with $p_g=q=1$ is base-point free, whence the bicanonical map $\phi:=\phi_{|2K|}: S \longrightarrow \mathbb{P}^{K_S^2}$ of $S$ is a morphism. Moreover such a morphism is generically finite by [Xi85a], so  $\phi(S)$ is a surface $\Sigma$. We will say that a surface $S$ contains a \emph{genus 2 pencil} if there is a morphism $f: S \longrightarrow B$, where $B$ is a smooth curve and the general fibre $\Phi$ of $f$ is a smooth curve of genus $2$. Notice that in this case the bicanonical map $\phi$ of $S$ is not birational, since $|2K_S|$ cuts out on $\Phi$ a subseries of the bicanonical series of $\Phi$ which is composed with the hyperelliptic involution. In this case we say that $S$ presents the \emph{standard case} for the non-birationality of the bicanonical map; otherwise, namely if $\phi$ is not birational but $S$ does not contain any genus $2$ pencils, we say that $S$ presents the \emph{non-standard} case. By the results of Bombieri (later improved by Reider, see [Bo73] and [Re88] ) it follows that, if $K_S^2 \geq 10$ and the bicanonical map is not birational, then $S$ contains a genus $2$ pencil. Hence there exist only finitely many families of surfaces of general type presenting the non-standard case, and one would like to classify all of them; however, this problem is still open, although many examples are known. In the paper [Xi90] G. Xiao gave two lists of possibilities for the bicanonical image of such a surface; later on several authors investigated their real occurrence. For more details about this argument, we refer the reader to the paper [Ci97]. \\
No examples of surfaces with $p_g=q=1$ and presenting the non-standard case were hitherto known; if $S$ is such a surface and $K_S^2 \geq 5$, then a result of Xiao ([see Xi90, Proposition 5]) implies that the degree of $\phi$ is either $2$ or $4$. In this work we describe the surfaces of general type with $p_g=q=1, \; K_S^2=8$ and such that the degree of $\phi$ is $2$. It will turn out that they belong to three distinct families, which provide as well the first-known examples of surfaces which such invariants. None of these surfaces contains a genus $2$ pencil, thus they are three substantially new pieces in the classification of surfaces presenting the non-standard case. \\
What we will show is that, as in the case $p_g=q=0$, the surfaces with $p_g=q=1, \; K_S^2=8$ and bicanonical map of degree $2$ are isogenous to a product. More precisely, the aim of this paper is to prove the following Theorem \ref{main theorem}.

\begin{theorem} \label{main theorem}
Let $S$ be a minimal surface of general type with $p_g=q=1, \; K_S^2=8$ and such that its bicanonical map has degree $2$. Then $S$ is a quotient of type $S=(C \times F)/G$, where $C, \;F$ are smooth curves and $G$ is a finite group acting faithfully on $C, \; F$ and freely on $C \times F$. Moreover $C$ is a curve of genus $3$ which is both hyperelliptic and bielliptic, $C/G$ is an elliptic curve $E$ isomorphic to the Albanese variety of $S$ and $F/G \cong \mathbb{P}^1$. The bicanonical map $\phi$ of $S$ factors through the involution $\sigma$ of $S$ induced by the involution $\tau \times id$ on $C \times F$, where $\tau$ is the hyperelliptic involution of $C$. The occurrences for $g(F)$ and $G$ are the following three:

\begin{itemize}   
\item[$I$.] $g(F)=3, \; G \cong \mathbb{Z}_2 \times \mathbb{Z}_2$;
\item[$II$.] $g(F)=4, \; G \cong S_3$;
\item[$III$.]$g(F)=5, \; G \cong D_4$.
\end{itemize}
The curve $F$ is hyperelliptic in case $I$, whereas it is not hyperelliptic in cases $II$ and $III$. \\
Surfaces of type $I, \; II, \; III$ do exist and they form three generically smooth, irreducible components $\mathcal{S}_{I}, \; \mathcal{S}_{II}, \; \mathcal{S}_{III}$ of the moduli space $\mathcal{M}$ of surfaces with $p_g=q=1, \; K_S^2=8$, whose respective dimensions are:
\begin{displaymath}
\textrm{dim }\mathcal{S}_{I}=5, \quad \textrm{dim }\mathcal{S}_{II}=4, \quad \textrm{dim }\mathcal{S}_{III}=4. 
\end{displaymath}
\end{theorem}
It is worth remarking that Theorem \ref{main theorem} still holds if one replaces the assumption that the bicanonical map $\phi$ has degree $2$ with the assumption that it factors through a rational involution (see Theorem \ref{grado2fattor}).\\ 

Let us briefly describe the plan of this article. \\
In Section \ref{degree 2} we give some generalities about surfaces with bicanonical map of degree $2$, and we introduce a very useful method in order to study them: the analysis of the \emph{bicanonical involution}, following [Xi90] and [CM02]. Besides, we present the two lists of Xiao [Xiao' s list $A$] and [Xiao's list $B$] that yield the possibilities for the bicanonical image when $S$ presents the non-standard case. \\ 
In Section \ref{p_g=q=1} we show that if $S$ is a minimal surface of general type with $p_g=q=1, \; K_S^2=8$ and bicanonical map of degree $2$, then $S$ is isogenous to a product, i.e. $S=(C \times F)/G$. We show moreover that there are at most three families of such surfaces, and we describe them. \\
In Section \ref{tech tool} we introduce some technical tools needed for what follows. This section is divided into two parts. In the former part we briefly study the group of automorphisms of a hyperelliptic curve, and we present a theorem of Accola (Theorem \ref{accola}); these topics will be useful in order to give an explicit equation for the curve $C$. In the latter part we deal with the main tools that will allow us to construct the curve $F$, and to show that it varies in an irreducible family; in particular we will need the orbifold fundamental group and the Hurwitz monodromy. \\  
In Section \ref{esempi} we show that the three families described in Section \ref{p_g=q=1} actually exist, by constructing the two curves $C, \; F$ and by exhibiting explicitly the actions of $G$ on them. \\ 
In Section \ref{bicanonica mappa} we prove that for any surface $S$ constructed in Section \ref{esempi}, the bicanonical map $\phi$ has degree $2$, and it is composed with the involution $\sigma$ of $S$ induced by $\tau \times id$.\\
In Section \ref{modellipiani} we give an alternative description of the three families as generically double coverings of $\mathbb{P}^2$, recovering a construction proposed by Du Val in [DV52]. \\
In Section \ref{mmoduli} we show that the surfaces with bicanonical map of degree $2$ form three disjoint, smooth, irreducible components $\mathcal{S}_I, \; \mathcal{S}_{II}, \; \mathcal{S}_{III}$ of the moduli space $\mathcal{M}$ of surfaces with $p_g=q=1, \; K_S^2=8$, whose dimensions are respectively $5, \; 4, \; 4$. \\
Finally, in Section \ref{open problems} we present some open problems. \\ \\
$\mathbf{Notations \; and \; conventions}$ All varieties, morphisms, etc. in this article are defined over the field $\mathbb{C}$ of the complex numbers. By ``surface'' we mean a projective, non-singular surface $S$, and for such a surface $K_S$ denotes the canonical class, $p_g(S)=h^0(S, \; K_S)$ is the \emph{geometric genus},  $q(S)=h^1(S, \; K_S)$ is the \emph{irregularity} and $\chi(\mathcal{O}_S)=1-q(S)+p_g(S)$ is the \emph{Euler characteristic}. We say that $S$ is \emph{regular} if $q(S)=0$, \emph{irregular} otherwise. We denote by $\phi$ the \emph{bicanonical map} of $S$, namely the rational map $\phi: S \longrightarrow \mathbb{P}^{K_S^2 +\chi(\mathcal{O}_S)-1}$ associated to the complete linear system $|2K_S|$. If $\sigma: S \longrightarrow S$ is an involution and $\psi: S  \longrightarrow S / \sigma$ is the quotient map, then we will say that $\phi$ is \emph{composed} with $\sigma$ if $\phi$ factors through the double cover $\psi$.   \\ \\     
$\mathbf{Acknowledgements.}$ It is a pleasure to thank C. Ciliberto
for his generous help and constant encouragement during the
preparation of this work. Part of this research was done during the
spring semester of the academic year 2002-2003, when the author was a
visiting student of the University of Bayreuth (Germany), and was
supported by Eu Research Training Network EAGER,
no. HPRN-CT-2000-00099. In particular he thanks F. Catanese and
I. Bauer for many useful discussions and suggestions and for their
supportive attitude. Also, he is indebted to F. Catanese, who
suggested the correct proof for the irreducibility of the moduli space
of surfaces of type $III$, and to R. Pardini, who read an earlier
version of this paper and pointed out some mistakes.
\section{Surfaces with bicanonical map of degree $2$} \label{degree 2}
\subsection{The bicanonical involution} \label{bicanonic inv}
Let $S$ be a minimal surface of general type and let $\sigma: S \longrightarrow S$ be a birational involution. Then $\sigma$ is biregular, and its fixed locus is given by a (possibly reducible) smooth curve $R'$ and isolated fixed points $p_1, \ldots, p_t$. Let $\pi: \hS \longrightarrow S$ be the blow-up of $S$ at $p_1, \ldots, p_t$; therefore $\sigma$ extends to an involution $\hat{\sigma}: \hS \longrightarrow \hS$ whose fixed locus is
\begin{displaymath}
\hR = \hR' + \sum_{i=1}^t E_i,
\end{displaymath}
where $E_i$ is the exceptional divisor over $p_i$ and $\hR'$ is a smooth curve isomorphic to $R'$. Let $T:=S / \sigma, \; \hW:=\hS / \hat{\sigma}$ and let 
$ \psi: S \longrightarrow T, \;  \hat{\psi}: \hS \longrightarrow \hW$ be the projections onto the quotients. The surface $T$ has $t$ nodes, whereas $\hW$ is smooth and we have a commutative diagram
\begin{displaymath} 
\begin{CD}
\hS  @>\pi>> S\\
@VV{\hat{\psi}}V  @VV{\psi}V\\
\hW @> \rho>> T,\\
\end{CD}
\end{displaymath}
where $\rho$ is the blow-up of $T$ at the nodes. $\hat{\psi}$ is a double cover and its branch locus $\hB$ is given by
\begin{displaymath}
\hB= \hB'+ \sum_{i=1}^t \Omega_i,
\end{displaymath}
where the $\Omega_i$'s  are $(-2)-$curves. Let $\hmL$ be the element in Pic$(\hW)$ such that $2 \hmL=\hB$ and which determines $\hat{\psi}$. Then we have $\hat{\psi}_* \mathcal{O}_{\widehat{S}}=\mathcal{O}_{\widehat{W}} \oplus \hmL^{-1}$, where $\mathcal{O}_{\widehat{S}}$ is the invariant part and $\hmL^{-1}$ is the anti-invariant part of $\hat{\psi}_*\mathcal{O}_{\widehat{S}}$ under the action of $\hat{\sigma}$. Since $\hat{\psi}$ is a double cover, the invariants of $\widehat{S}$ and $\widehat{W}$ are related in the following way (see [BPV84, p.183]):
\begin{displaymath} 
\begin{array}{ll}
 K_{\widehat{S}}^2 & = \quad 2(K_{\widehat{W}}+\hmL)^2,\\
 \chi(\mathcal{O}_{\widehat{S}}) & = \quad 2 \chi(\mathcal{O}_{\widehat{W}}) + \frac{1}{2} \hmL \cdot (K_{\widehat{W}}+\hmL),\\
 p_g(\widehat{S}) & = \quad p_g(\widehat{W})+h^0(\widehat{W}, \; K_{\widehat{W}}+\hmL),\\
 q(\widehat{S}) & = \quad q(\widehat{W})+h^1(\widehat{W}, \; K_{\widehat{W}}+\hmL).
\end{array}
\end{displaymath}
Let us now consider the bicanonical map $\phi$ of $S$; the following result is due to Mendes Lopes and Pardini (see [MP03, Proposition 2.1]).
\begin{proposition} \label{fattorizzazione}
Let $S$ be a minimal surface of general type such that
$\chi(\mathcal{O}_S)=1$, and let $\sigma$ be an involution of $S$. Let
$T:=S/ \sigma$ be the quotient, and let $\psi: S \longrightarrow T$ be the projection. Suppose moreover that $|2K_S|$ has no fixed components, and that $p_g(T)=q(T)=0$. Then $\phi$ is composed with $\sigma$ if and only if the number $t$ of isolated fixed points of $\sigma$ is equal to $K_S^2 +4$.
\end{proposition}
\begin{proof}
Let $\pi: \hS \longrightarrow S$ be the blow-up of $S$ at the isolated
fixed points of $\sigma$, and let $\hat{\sigma}$ be the involution induced by $\sigma$ on $\hS$. Using the previous notations we can write
\begin{displaymath} 
H^0(\hS, \; 2K_{\hS})= H^0( \hW, \; 2K_{\hW} + \widehat{\mathcal{L}}) \oplus H^0(\hW, \; 2K_{\hW}+ \hB'+ \sum_{i=1}^t \Omega_i),
\end{displaymath} 
where $\hat{\psi}^* H^0( \hW, \; 2K_{\hW} + \widehat{\mathcal{L}})$ is the $(-1)-$eigenspace for the action of $\hat{\sigma}$ on $H^0(\hS, \; 2K_{\hS})$, whereas $\hat{\psi}^* H^0(\hW, \; 2K_{\hW}+ \hB'+ \sum_{i=1}^t \Omega_i)$ is the $(+1)-$eigenspace. The map $\phi$ factors through $\sigma$ if and only if one of these eigenspaces is equal to zero. On the other hand, we have
\begin{displaymath}
|2K_{\hS}|=\pi^*|2K_S|+2 \sum_{i=1}^t E_i,
\end{displaymath}
and this in turn means that the fixed part of $|2K_{\hS}|$ is given by $\sum_{i=1}^t E_i$, since we are supposing that $|2K_S|$ has no fixed part. It follows that there exists some bicanonical curve on $\hS$ which does not contain the curve $\hR'$. Observe now that all the ``anti-invariant'' sections of $|2K_{\hS}|$ must contain $\hR'$; this implies that $\phi$ factors through $\sigma$ if and only if $H^0( \hW, \; 2K_{\hW} + \widehat{\mathcal{L}})=0$, i.e. if and only if $H^0(\hS, \; 2K_{\hS})=H^0(\hW, \; 2K_{\hW}+ \hB'+ \sum_{i=1}^t \Omega_i)$. Now we have to consider two cases.\\ \\
$- \; \;  \hB'=0$. \\
In this case $K_S= \psi^*K_T$, hence $K_T$ is nef and big because $K_S$ is nef and big by assumption, and this implies that $T$ is a surface of general type. The map $\hat{\psi}: \hS \longrightarrow \hW$ is branched only along $t$ $(-2)-$curves, and the standard formulae for double covers give
\begin{displaymath}
1= \chi(\mathcal{O}_S)=2 \chi(\mathcal{O}_T) - \frac{t}{4}=2-\frac{t}{4},
\end{displaymath} 
that is $t=4 < K_S^2+4$. Moreover
\begin{equation} \label{pidue}  
P_2(S)=K_S^2+1=2K_T^2+1 >K_T^2+1=P_2(T).
\end{equation}
On the other hand, if $\phi$ factors through $\sigma$ then we get $H^0(\hS, \; 2K_{\hS})=H^0(\hW, \; 2K_{\hW}+ \sum_{i=1}^t \Omega_i)$, and this in turn implies $H^0(S, \; 2K_S)= H^0(T, \; 2K_T)$, that is, $P_2(S)=P_2(T)$, which contradicts (\ref{pidue}). Hence this case does not occur. \\ \\
$- \; \; \hB' \neq 0$. \\
Recall that $2K_S= \psi^*(2K_T+B)$, hence $2K_T+B$ is nef and big. We have the following equality of $\mathbb{Q}-$divisors:
\begin{displaymath}  
K_{\hW}+ \widehat{\mathcal{L}}=\frac{1}{2}(2K_{\hW}+\hB')+ \frac{1}{2} \sum_{i=1}^t \Omega_i.
\end{displaymath}
Notice that the divisor $\frac{1}{2}(2K_{\hW} + \hB')= \frac{1}{2} \psi^*(2K_T+B)$ is nef and big, whereas the $\mathbb{Q}-$divisor $\frac{1}{2} \sum_{i=1}^t \Omega_i$ is effective with normal crossing support and zero integral part. Thus $h^i( \hW, \; 2K_{\hW} + \widehat{\mathcal{L}})=0$ for $i>0$ by the Kawamata-Viehweg theorem, and so
\begin{equation} \label{unorel}
\begin{split}
h^0(\hW, \; 2K_{\hW}+ \widehat{\mathcal{L}}) &= \chi(2K_{\hW}+ \widehat{\mathcal{L}}) \\
&=1+(K_{\hW})^2+ \frac{3}{2}K_{\hW} \widehat{\mathcal{L}}+ \frac{1}{2} (\widehat{\mathcal{L}})^2.
\end{split} 
\end{equation} 
On the other hand, projection formula for double covers yields
\begin{equation} \label{duerel}
1= \chi(\mathcal{O}_{\hS})= \chi(\mathcal{O}_{\hW})+ \chi(K_{\hW}+ \widehat{\mathcal{L}})= 1 + \chi(K_{\hW}+ \widehat{\mathcal{L}}),
\end{equation} 
which by Riemann-Roch is equivalent to
\begin{equation} \label{trerel}
(\widehat{\mathcal{L}})^2+ K_{\hW} \widehat{\mathcal{L}}=-2.
\end{equation}
Finally we have
\begin{equation} \label{quatrel}
t=K_S^2-(K_{\hS})^2=K_S^2-2(K_{\hW}+ \widehat{\mathcal{L}})^2.
\end{equation}  
Using equalities (\ref{unorel}), (\ref{duerel}) and (\ref{trerel}), relation (\ref{quatrel}) can be rewritten as
\begin{displaymath}
t=K_S^2+4-2h^0(2K_{\hW}+ \widehat{\mathcal{L}}).
\end{displaymath}
Summing up, $\phi$ is composed with $\sigma$ if and only if $h^0(2K_{\hW}+ \widehat{\mathcal{L}})=0$, i.e. if and only if $t=K_S^2+4$.
\end{proof}
Now suppose that $S$ is a minimal surface of general type such that $|2K_S|$ is without fixed components and the bicanonical map $\phi$ of $S$ has degree $2$. Then $S$ possess the so-called \emph{bicanonical involution} $\sigma: S \longrightarrow S$, which exchanges the two sheets of the (generically) double cover $\phi: S \longrightarrow \Sigma$. Notice that by definition $\phi$ is composed with $\sigma$, hence we have a commutative diagram
\begin{equation} \label{triang 5}
\begin{CD} 
\xymatrix{
\hS \ar[rrrr]^{\hat{\phi}} \ar[rrd]_{\hat{\psi}} & & & & \Sigma \\
& & \widehat{W} \ar[rru]_{\mu} & &
} 
\end{CD}
\end{equation}
where the map $\mu$ is birational. 
The existence of such a diagram implies that the following equalities hold (see [CM02, Proposition 6.1]):
\begin{equation} \label{relazionirivestimento}
\begin{array}{ll}
(i) & (2K_{\widehat{W}}+\hB')^2=2K_S^2;\\
(ii) & \chi(\mathcal{O}_{\widehat{W}}(2K_{\widehat{W}}+\hmL))=0;\\
(iii) & K_{\widehat{W}} \cdot (K_{\widehat{W}}+\hmL)= \chi(\mathcal{O}_{\widehat{W}})-\chi(\mathcal{O}_S).
\end{array}
\end{equation}
\begin{proposition} \label{relazioni}
Let $S$ be such that the degree of $\phi$ is $2$, let $\sigma$ be the
bicanonical involution and let $\hW=\hS / \hat{\sigma}$. Suppose moreover $\chi(\mathcal{O}_{\hW})= \chi( \mathcal{O}_S)=1$. Then we have:
\begin{itemize}
\item[($a$)] $t=K_S^2+4$;
\item[($b$)] $K_{\hW}\hB'=K_{\hW}\hB=-2K_{\hW}^2$;
\item[($c$)] $(\hB')^2=4K_{\hW}^2+2K_S^2$;
\item[($d$)]$(R')^2=2K_{\hW}^2+K_S^2$;
\item[($e$)] $K_SR'=K_S^2$. 
\end{itemize}
\end{proposition}
\begin{proof}
Relation ($a$) comes from Proposition \ref{fattorizzazione}.
Since $2 \hmL=\hB'+ \sum_{i=1}^t \Omega_i$, we apply formulae (\ref{relazionirivestimento}) obtaining
\begin{displaymath}
\begin{array}{ll}
(i) & 4K_{\hW}^2+4K_{\hW}\hB'+(\hB')^2=2K_S^2;\\
(ii) & 8K_{\hW}^2+6K_{\hW}\hB'+(\hB')^2=2t-8;\\
(iii) & 2K_{\hW}^2+K_{\hW} \hB'=0.
\end{array}
\end{displaymath}
The last equation implies $K_{\hW}\hB'=-2K_{\hW}^2$ which is ($b$); substituting in the first two relations we obtain
\begin{displaymath} 
\begin{array}{ll}
(i') & -4K_{\hW}^2+(\hB')^2=2K_S^2;\\
(ii') & -4K_{\hW}^2+(\hB')^2=2t-8.
\end{array}
\end{displaymath}
Relation $(i')$ yields $(\hB')^2=4K_{\hW}^2+2K_S^2$, and in this way we proved ($c$).  Finally we prove ($d$) and ($e$). Note that, since $\hB'$ and $\hR'$ are isomorphic curves, we have $p_a(\hR')=p_a(\hB')$, and the right-hand side is equal to $K_{\hW}^2+K_S^2+1$ by ($b$) and ($c$). Therefore we have $(\hR')^2=\frac{1}{2}(\hB')^2=2K^2_{\hW}+K_S^2$, which is ($d$), and we can calculate $K_SR'$ obtaining $K_SR'=K_S^2$, which is ($e$).
\end{proof}  

\subsection{The good minimal model} \label{sec: HX}
Let $S$ be a minimal surface of general type such that its bicanonical map has degree $2$, and let $\sigma$ and $\hW$ be as in Proposition \ref{relazioni}. Assume moreover that $\hW$ is a \emph{ruled} surface.

\begin{definition} \label{HX-minimal}
A \emph{good minimal model} of $\hW$ is a minimal model $p: \hW \longrightarrow W$ with a smooth ruling $\theta: W \longrightarrow C$ such that the following properties are verified:  
\begin{enumerate}
\item Let $B$ be the image of $\hB$ on $W$. Then the degree $k$ of $B$ over $C$ is minimal along all such choices.
\item The greatest order of singularities of $B$ is minimal, and the number of singularities of $B$ with greatest order is minimal, among all the choices satisfying condition $(1)$.
\end{enumerate} 
\end{definition}

\begin{proposition} \label{esiste HX}
If $\hW$ is not isomorphic to $\mathbb{P}^2$, then $\hW$ admits a good minimal model $W$.
\end{proposition}
\begin{proof}
The hypothesis that $\hW$ is not isomorphic to $\mathbb{P}^2$ implies that $\hW$ admits a birational morphism $p_0: \hW \longrightarrow W_0$, where $W_0$ is geometrically ruled. We consider all the products $g$ of elementary transformations such that $g \circ p_0$ is a morphism. Among all such $g$, we chose one that first minimizes $k$, and second minimizes the order and the number of singularities of $B$. Then $p:= g \circ p_0: \hW \longrightarrow W$ is a good minimal model of $\hW$.
\end{proof}
On the other hand, the case $\hW=\mathbb{P}^2$ is very easy to describe.

\begin{proposition} \label{W=P2}
If $\hW$ is isomorphic to $\mathbb{P}^2$, then we have the following two possibilities:
\begin{itemize}
\item $p_g=6, \thinspace q=0, \thinspace K_S^2=8$.\\
In this case $S$ is the double cover of the plane $\psi: S \longrightarrow \mathbb{P}^2$, branched along a smooth curve $B$ of degree $10$; 
\item $p_g=3, \thinspace q=0, \thinspace K_S^2=2$.\\
In this case $S$ is the double cover of the plane $\psi: S \longrightarrow \mathbb{P}^2$, branched along a smooth curve $B$ of degree $8$.
\end{itemize}
\end{proposition}  
\begin{proof}
Since $\mathbb{P}^2$ does not contain $(-2)-$curves, we have $t=0$, so $S=\hS$  and there is a double cover $\psi: S \longrightarrow \mathbb{P}^2$. Moreover, $S$ is smooth, so the branch locus $B$ of $\psi$ is a smooth divisor; it follows in particular that $S$ is a regular surface. Let $2n$ be the degree of the branch curve $B \subset \mathbb{P}^2$; the fact that the bicanonical map of $S$ factors through $\psi$ implies $\chi(\mO_{\mathbb{P}^2}(2K_{\mathbb{P}^2}+nH))=0$, where $H$ is the class of a line. This gives $n^2-9n+20=0$, therefore either $n=4$ or $n=5$, and this in turn implies that $S$ is the double cover of the projective plane branched along a smooth divisor $B$ whose degree is either $8$ or $10$.
\end{proof}

\begin{remark} 
The idea of good minimal model goes back to Xiao, see $[$Xi90, p.725$\;]$. However a quite similar definition was given by Hartshorne in $[$Ha69, Proposition 3.1$\;]$ in order to study curves with high self-intersection on a ruled surface.
\end{remark}

\begin{remark}
The examples in Proposition \ref{W=P2} were first studied by Du Val in $[$DV52$\;]$. In this paper he investigated regular surfaces of general type which present the non-standard case, under the hypothesis that the general canonical curve is smooth and irreducible. For a discussion about Du Val's results see $[$Ci97$\;]$.  
\end{remark}  
Proposition \ref{W=P2} shows that we can suppose without loss of
generality that $\hW$ is not isomorphic to $\mathbb{P}^2$; then
Proposition \ref{esiste HX} says that we can consider a good minimal
model $W$ of $\hW$.  $W$ is a geometrically ruled surface over a curve
$C$ such that $g(C)=q(\hW)$, and the map $p:\hW \longrightarrow W$ is
composed of a series of blowing-ups. Let $x_1,x_2, \ldots, x_s$ be the
centers of these blowing-ups, and let $\mathcal{E}_i$ be the inverse image of $x_i$ on $\hW$ (with right multiplicities such that $\mathcal{E}_i \mathcal{E}_j=- \delta_{ij}, \; K_{\hW}=p^* K_W+\sum_{i=1}^s \mathcal{E}_i$). Then, according to [Xi90, p.725], we can write
\begin{equation} \label{hat B}
\hB=p^*B- \sum_{i=1}^sb_i \mathcal{E}_i,
\end{equation}
where the $b_i$'s are even positive integers.
\begin{definition} \label{singularity}
\begin{itemize}
\item[]
\item A \emph{negligible singularity} is a point $x_j$ such that $b_j=2$, and $b_i \leq 2$ for any point $x_i$ infinitely near to $x_j$.
\item A $[2b+1,2b+1]-$ \emph{singularity} is a pair $(x_i,x_j)$ such that $x_i$ is immediately infinitely near to $x_j$ and $b_i=2b+2, \thinspace b_j=2b$.
\end{itemize}
\end{definition}
For example, a double point or an ordinary triple point are negligible singularities, whereas a $[3,3]-$point is not.  
Moreover, we will say that a $[2b+1,2b+1]-$singularity $\xi=(x_i,x_j)$ is \emph{tangent} to a fibre $L$ of $W$ if $x_j \in L$ and $x_i$ is infinitely near to $x_j$ in the direction given by $L$.\\
Let $W$ be a good minimal model of $\hW$ and let $L$ be a general fibre of $W$; if $k$ is the degree of $B$ over the base of the ruling and $C_0$ is the section of minimal self-intersection, then there exists an unique integer $l$ such that
\begin{equation} \label{espr B}
B \cong kC_0+ \Big( \frac{1}{2}ke+l \Big) L,
\end{equation}
where $e=-C_0^2$.
\begin{proposition} \label{relriv 5}
If $\chi(\mO_{\hW})= \chi(\mO_S)=1$, then the following relations hold:
\begin{enumerate} 
\item  $2kl- \sum_{i=1}^s b_i^2=4K_{\hW}^2-8$;
\item  $-2(k+l)+\sum_{i=1}^sb_i=-2K_{\hW}^2$.
\end{enumerate}
\end{proposition}
\begin{proof}
From Proposition \ref{relazioni} we obtain $K_{\hW}\hB=-2K_{\hW}^2$, $(\hB)^2=(\hB')^2-2t=4K_{\hW}^2-8$; hence it is equivalent to prove the following two relations:
\begin{enumerate}
\item $(\hB)^2 = 2kl- \sum_{i=1}^s b_i^2$;
\item $K_{\hW}\hB =-2(k+l)+ \sum_{i=1}^s b_i$.
\end{enumerate} 
Actually we have:
\begin{displaymath}
\begin{split}
(\hB)^2 & =\Big( p^*B -\sum_{i=1}^s b_i \mathcal{E}_i \Big) ^2\\
        & =B^2-\sum_{i=1}^s b_i^2= \Big( kC_0+ \Big( \frac{1}{2}ke+l \Big) L \Big)^2- \sum_{i=1}^s b_i^2\\
        & =2kl- \sum_{i=1}^s b_i^2,
\end{split}
\end{displaymath}
which is ($1$), and:
\begin{displaymath}
\begin{split}
K_{\hW}\hB &=\Big( p^*K_W+ \sum_{i=1}^s \mathcal{E}_i \Big) \Big( p^*B-\sum_{i=1}^s b_i \mathcal{E}_i \Big)\\
            &=K_WB+\sum_{i=1}^s b_i=\Big( -2C_0-(e+2)L \Big) \Big( kC_0+ \Big( \frac{1}{2}ke+l \Big) L \Big)+\sum_{i=1}^s b_i\\
            &=-2(l+k)+\sum_{i=1}^s b_i,
\end{split}
\end{displaymath}
which is ($2$). This completes the proof. 
\end{proof}

\subsection{The two lists of Xiao}
The paper [Xi90] is devoted to classify all possible non-standard cases for the non-birationality of the bicanonical map of a minimal surface of general type $S$. Xiao gave a long list of possibilities; later on several authors investigated their real occurrence.\\
Here we suppose that $S$ is \emph{irregular} and such that its bicanonical map has degree $2$. We first consider the case when $\hW$ is ruled.  In [Xi90, Theorem 2] and [Xi90, Proposition 6] Xiao proved the following result:

\begin{theorem}[$\mathbf{Xiao's \; list \; A}$] \label{Xiao list 1}
Let $S$ be a minimal irregular surface of general type with
bicanonical map of degree $2$ and presenting the non-standard
case. Suppose moreover that $\hW$ is a ruled surface. Then $\hW$ is
$\emph{rational}$. Let $W$ be a good minimal model of $\hW$, and let $k, \;l$ be as in the previous section. Then only the following possibilities can occur:
\begin{itemize}
\item[($A1$)] $k=16, \thinspace l=18$.\\  
$B$ contains three singularities $[9,9]$, a singularity of order $8$ and no other singularities of order greater than $4$. In this case $S$ contains a pencil of hyperelliptic curves of genus $7$ with three double fibres.
\item[($A2$)] $k=12, \thinspace  l=14$.\\
$B$ contains three singularities $[7,7]$ and no other singularities of order greater than $4$. In this case $S$ contains a pencil of hyperelliptic curves of genus $5$ with three double fibres.
\item[($A3$)] $k=8, \thinspace l=8+2i$, where $1 \leq i \leq 5$.\\
$B$ contains $i+1$ singularities $[5,5]$, plus possibly other singularities of strictly smaller order. In this case $S$ contains a pencil of hyperelliptic curves of genus $3$ with $i+1$ double fibres.
\item[($A4$)] $k=8, \; l=6$.\\
If this case occurs, then $S$ verifies $p_g=q=1, \; K_S^2=3$. $B$ contains six singularities $[3,3]$, plus possibly other negligible singularities, and $S$ contains a pencil of hyperelliptic curves of genus $3$ with six double fibres.
\end{itemize}
\end{theorem}      

\begin{proposition} \label{fibre in ram}
Let $S$ be as in the previous Theorem \ref{Xiao list 1}, and let
$\xi=(x_i,x_j)$ be one of the $[\frac{1}{2}k+1, \frac{1}{2}k+1]-$singularities of $B$. Let $L_{\xi}$ be the fibre of $W$ containing $\xi$. Then $\xi$ is tangent to $L_{\xi}$ and $L_{\xi} \subset B$. As a consequence, $B$ does not contain $C_0$.
\end{proposition} 
\begin{proof}
Suppose that the singularity $\xi$ is not tangent to $L_{\xi}$. This means that $B$ does not contain $L_{\xi}$ and that the intersection multiplicity of $B$ and $L_{\xi}$ in $x_j$ is exactly $\frac{1}{2}k+1$. Since $BL=k$, after performing an elementary transformation centered at $x_j$ we obtain another minimal model $W'$, where the singularity $\xi$ is substituted by a singular point of order $\frac{1}{2}k-1$. This contradicts the fact that $W$ is a good minimal model, so $\xi$ must be tangent to $L_{\xi}$. Then the intersection multiplicity of $B$ and $L$ at $x_j$ verifies $\textrm{mult}_{x_j}(B,L)=k+2>k$, and this in turn implies $L_{\xi} \subset B$.
\end{proof} 

\begin{corollary} \label{about C0}
We can assume that the curve $C_0$ does not meet the $[\frac{1}{2}k+1, \frac{1}{2}k+1]-$singularities of $B$.
\end{corollary}
\begin{proof}
Let $\xi:=(x_i,x_j)$ be one of the $[\frac{1}{2}k+1,\frac{1}{2}k+1]-$singularities of $B$, and suppose that $C_0$ contains $\xi$. Of course the point $x_j$ is a fundamental point of the map $p: \hW \longrightarrow W$, hence if  $\textrm{elm}_{x_j}: W \dashrightarrow W'$ is the elementary transformation centered at $x_j$, the birational map $\textrm{elm}_{x_j} \circ p$ is a morphism. In this way we did not change the value of $k$, and we found another good minimal model $W'$ with the property that the number of $[\frac{1}{2}k+1, \frac{1}{2}k+1]-$singularities of $B$ contained in $C_0$ is dropped by $1$. Now we can repeat this process and after a finite number of steps we are done. 
\end{proof}

\begin{corollary} \label{fibre doppie}
Every $[ \frac{1}{2}k+1, \frac{1}{2}k+1]-$singularity of $B$ gives rise to two disjoint $(-1)-$curves in $\widehat{S}$, both contained in the fixed locus of $\hat{\sigma}$.
\end{corollary}
\begin{proof}
This follows from the \emph{canonical resolution} of the singularities of a double cover; see [BPV84, p.87].
\end{proof}
Now we consider the case when  $\hW$ is not ruled. If this happens, [Xi90, Theorem $3$], [CCM98, Theorem $A$] and [CM02, Theorem $1.1$] give the following list:

\begin{theorem}[$\mathbf{Xiao's \; list\; B}$] \label{Xiao li 1,2}
Let $S$ be a minimal irregular surface of general type, with bicanonical map of degree $2$ and presenting the non-standard case. Suppose that $\hW$ is not a ruled surface, and let $W$ be the $($unique$)$ minimal model of $\hW$. Then only the following possibilities can occur:
\begin{itemize}
\item[($B1$)] $\chi(\mO_S)=1$, $K_S^2=3$ or $4$, $W$ is an Enriques surface.
\item[($B2$)] $p_g=q=1$, $3 \leq K_S^2 \leq 6$, $W$ is a regular surface with geometric genus $1$, and with Kodaira dimension $1$.
\item[($B3$)] $p_g=q=2$, $K_S^2=4$.
\item[($B4$)] $p_g=q=3$, $K_S^2=6$.
\end{itemize}
\end{theorem}
We remark that cases ($B3$) and ($B4$) in the above list really occur. In fact in [CM02, Theorem 1.1] it is shown that a surface of type ($B3$) is a double cover of a principally polarized abelian surface $(A, \Theta)$, with $\Theta$ irreducible, branched along a divisor $B \in |2 \Theta|$ having at most rational double points as singularities, whereas in [CCM98, Theorem $A$] it is shown that a surface of type ($B4$) is the double symmetric product of a smooth curve of genus $3$.  \\ \\
Now we give some computation in the case $\chi(\mO_S)=1$ that will be useful in what follows:

\begin{proposition} \label{Xiao list 2}
Let $S$ be a minimal surface of general type with $\chi(\mO_S)=1$ and bicanonical map $\phi$ of degree $2$, presenting the nonstandard case. Suppose moreover that $\hW$ is ruled, in such a way that Theorem \ref{Xiao list 1} applies (in particular $\hW$ is rational). Let $\delta$ be the number of $b_i$ which are equal to $2$ in $($\ref{hat B}$)$. Then the following holds:
\begin{itemize}
\item if $S$ belongs either to case $(A1)$ or to case $(A2)$ of $[$Xiao's list $A]$, then we have
\begin{center}
$(R')^2=K_S^2-2 \delta-2$;
\end{center}  
\item if $S$ belongs to case $(A3)$ of $[$Xiao's list $A]$, then we have
\begin{center}
$(R')^2=K_S^2-2 \delta -2i +2$;
\end{center}
\item if $S$ belongs to case $(A4)$ of $[$Xiao's list $A]$, then we have
\begin{center}
$(R')^2=K_S^2 - 2 \delta-8$.
\end{center}
\end{itemize}
\end{proposition}
\begin{proof}
Since $\hW$ is rational, we have $\chi(\mathcal{O}_{\hW})=\chi(\mathcal{O}_S)=1$, so Proposition \ref{relazioni} and Proposition \ref{relriv 5} apply. \\ \\
$-$ Suppose that $S$ belongs to case ($A1$) of [Xiao's list $A$]. Then we have $k=16, \thinspace l=18$. Denote by $x_1, \ldots, x_s$ the nonnegligible singularities of $B$. Without loss of generality, we can suppose that $(x_1,x_2),(x_3,x_4), (x_5,x_6)$ are the three singularities of type $[9,9]$ and that $x_7$ is the singularity of order $8$. Therefore the relations in Proposition \ref{relriv 5} become
\begin{equation} \label{frelriv 6}
\left\{ \begin{array}{ll}
28- \sum_{j=8}^s b_j^2-4 \delta= 4K_{\hW}^2\\
-6+ \sum_{j=8}^s b_j +2 \delta = -2 K_{\hW}^2,
\end{array} \right.
\end{equation}
where $b_j \geq 4$ for any $j$, or $b_j=0$ for any $j$.
Since every $b_j$ is an even positive integer, there exists $\beta_j \in \mathbb{N}$ such that $b_j=2 \beta_j$. Substituting in (\ref{frelriv 6}) we obtain
\begin{equation} \label{frelriv 7}
\left\{ \begin{array}{ll}
7 - \sum_{j=8}^s \beta_j^2 - \delta= K_{\hW}^2 \\
3- \sum_{j=8}^s \beta_j - \delta = K_{\hW}^2.
\end{array} \right.
\end{equation}
This implies
\begin{displaymath}
\sum_{j=8}^s \beta_j(\beta_j-1)=4,
\end{displaymath}
and the only solution is $s=9, \thinspace \beta_8= \beta_9 =2$, hence $b_8=b_9=4$. Moreover, the second equation in (\ref{frelriv 7}) gives
\begin{displaymath}
K_{\hW}^2=3-\beta_8 -\beta_9- \delta= -1- \delta,
\end{displaymath}
hence $(R')^2=K_S^2- 2 \delta -2$ by Proposition \ref{relazioni}. \\ \\
$-$ Suppose that $S$ belongs to case $(A2)$ of [Xiao's list $A$]. Then $k=12, \thinspace l=14$. Let $(x_1,x_2), (x_3,x_4), (x_5,x_6)$ be the three $[7,7]-$singularities of $B$. We obtain in this case
\begin{displaymath} 
\left\{ \begin{array}{ll}
11- \sum_{j=7}^s \beta_j^2 - \delta = K_{\hW}^2 \\
5 - \sum_{j=7}^s \beta_j - \delta =K_{\hW}^2.
\end{array} \right.
\end{displaymath}
This implies
\begin{displaymath}
\sum_{j=7}^s \beta_j( \beta_j-1)=6,
\end{displaymath}
and since $B$ does not contain singularities of order greater than $4$ besides its three points $[7,7]$, the only possibility is $s=9, \thinspace  \beta_7= \beta_8 = \beta_9=2$. Moreover
\begin{displaymath}
K_{\hW}^2=5-\beta_7 -\beta_8 -\beta_9 - \delta=-1- \delta,
\end{displaymath}
so Proposition \ref{relazioni} gives $(R')^2=K_S^2 -2 \delta -2$ in this case. \\ \\
$-$ Suppose that $S$ belongs to case $(A3)$ of [Xiao's list $A$]. Then we have $k=8, \thinspace l=8+2i$ for some $1 \leq i \leq 5$. Let $(x_1,x_2), (x_3,x_4), \ldots, (x_{2i+1}, x_{2i+2})$ be the $i+1$ points of type $[5,5]$ of $B$. Then we have
\begin{displaymath} 
\left\{ \begin{array}{ll} 
21-5i- \sum_{j=i+2}^s \beta_j^2 - \delta = K_{\hW}^2 \\
11-3i -\sum_{j=i+2}^s \beta_j - \delta = K_{\hW}^2. 
\end{array} \right.
\end{displaymath} 
This implies
\begin{displaymath}
\sum_{j=i+2}^s \beta_j (\beta_j -1)= 10-2i,
\end{displaymath}
and since $B$ does not contain singularities of order greater than $4$ besides its points $[5,5]$, the only possibility is $s=6, \thinspace \beta_{i+2}= \ldots =\beta_6=2$. Moreover
\begin{displaymath}
K_{\hW}^2=11-3i-2(5-i)- \delta= 1-i- \delta,
\end{displaymath}
so Proposition \ref{relazioni} gives $(R')^2=K_S^2- 2 \delta -2 i +2$ in this case. \\ \\
$-$ Finally suppose that $S$ belongs to case $(A4)$ of [Xiao's list $A$]. Then we have $k=8, \; l=6$. Let $(x_1,x_2), \;(x_3,x_4), \ldots,(x_{11},x_{12})$ be the six points $[3,3]$ of $B$. We obtain
\begin{displaymath} 
\left\{ \begin{array}{ll} 
-4- \sum_{j=7}^s \beta_j^2 - \delta = K_{\hW}^2 \\
-4- \sum_{j=7}^s \beta_j - \delta = K_{\hW}^2. 
\end{array} \right.
\end{displaymath} 
This implies
\begin{displaymath}
\sum_{j=7}^s \beta_j (\beta_j -1)= 0.
\end{displaymath}
If $\beta_j \geq2$ for some $j$ we obtain a contradiction, then we must have $\beta_j=0$ for any $j$, according with the fact that $B$ does not contain nonnegligible singularities apart from the six points $[3,3]$. Then $K_{\hW}^2=- \delta -4$, hence Proposition \ref{relazioni} gives $(R')^2=K_S^2- 2\delta-8$ in this case.
\end{proof}

\begin{remark}
G. Borrelli communicated to us that he was recently able to rule out cases $(A1)$ and $(A2)$ in $[$Xiao list $A]$. See $[Bor03]$ for further details.
\end{remark}

\begin{remark} \label{vale rat inv}
A careful analysis of the Xiao's proof shows that Theorem \ref{Xiao list 1} and Proposition \ref{Xiao list 2} still hold if one replaces the assumption that the bicanonical map $\phi$ has degree $2$ with the assumption that $\phi$ factors through a rational involution.
\end{remark}

\section{The case $p_g=q=1$, $K_S^2=8$} \label{p_g=q=1}

\subsection{The isotrivial pencil}

Let $S$ be a minimal surface of general type with $p_g=q=1$, and let $\phi: S \longrightarrow \mathbb{P}^{K_S^2}$ be its bicanonical map.
\begin{proposition} \label{K_S^2 geq 5}
If $p_g=q=1$, $K_S^2 \geq 5$ and $S$ presents the non-standard case, then the degree of $\phi$ is either $2$ or $4$. Moreover if $K_S^2=9$, then the non-standard case does not occur.
\end{proposition}
\begin{proof}
The former statement is [Xi90, Proposition 5]; the latter one is [CM02, Proposition 3.1].  
\end{proof}

\begin{proposition} \label{caso standard}
\begin{enumerate}
\item[]
\item If $p_g=q=1$ and $S$ contains a \emph{rational} genus $2$ pencil, then either $K_S^2=2$ or $K_S^2=3$.
\item If $p_g=q=1$ and the Albanese pencil has genus $2$, then $2 \leq K_S^2 \leq 6$.
\end{enumerate}
\end{proposition}
\begin{proof}
See [Xi85b, Corollary 3 p.51 and Theorem 2.2 p.17].
\end{proof}
On the other hand, since $q=1$, the Albanese pencil is the unique irrational pencil on $S$; therefore we have:

\begin{corollary} \label{nonstd}
If $p_g=q=1, \; K_S^2=7,8$ and the bicanonical map of $S$ is not birational, then $S$ presents the non-standard case.
\end{corollary} 

\begin{proposition} \label{kappa ampio}
If $p_g=q=1$ and $K_S^2=8, \; 9$ then $K_S$ is ample.
\end{proposition}
\begin{proof}
This follows form Miyaoka's inequality, see [Mi84, Section 2].
\end{proof}

\begin{lemma} \label{ramif a}
Let $S$ be a minimal surface of general type with $p_g=q=1$ and bicanonical map $\phi$ of degree $2$. Suppose that $S$ presents the non-standard case and that $\hW$ is ruled. Then $(R')^2 \leq 0$. Moreover if $(R')^2=0$, then there exist  $F_1, \ldots , F_j$ fibres in the Albanese pencil such that the support $(F_i)_{red}$ of $F_i$ is a smooth curve and
\begin{displaymath}
R'= \sum_{i=1}^j (F_i)_{red}.
\end{displaymath}
\end{lemma}
\begin{proof}
First remark that $\hW$ is rational because of Theorem \ref{Xiao list 1}. Let $\alpha: S \longrightarrow E$ be the Albanese pencil of $S$. The bicanonical involution $\sigma: S \longrightarrow S$ induces an involution $\sigma_E :E \longrightarrow E$. If the quotient $E / \sigma_E$ were again an elliptic curve, then $\hW$ would be a connected fibration with elliptic base, a contradiction. Then $E / \sigma_E \cong \mathbb{P}^1$, hence $\sigma_E(x)= -x$, and this implies that $R'$ is a disjoint union of components of fibres of $\alpha$, so the conclusion follows from Zariski's lemma (see [BPV84, p.90]) and from the smoothness of $R'$.
\end{proof}
We are now ready to show the following result, that is the first step in order to prove Theorem \ref{main theorem}.

\begin{proposition} \label{isotrivial pencil}
Let $S$ be a minimal surface of general type with $p_g=q=1$ and $K_S^2=8$ such that the bicanonical map has degree $2$. Then $S$ contains an isotrivial pencil $|C|$ of hyperelliptic curves of genus $3$ with six double fibres.
\end{proposition}
\begin{proof}
We first show that $S$ contains a base-point free pencil $|C|$ of curves of genus $3$. Proposition \ref{kappa ampio} shows that $K_S$ is ample, then $S$ does not contain $(-2)-$ curves; this means that the curve $B$ does not contain negligible singularities. Hence, using the same notations as in Proposition \ref{Xiao list 2}, we have $\delta=0$. Since $K_S^2=8$, $S$ is not in [Xiao's list $B$], hence it belongs to [Xiao's list $A$]. But if either case ($A1$) or case ($A2$) occurs, then applying Proposition \ref{Xiao list 2} we obtain $(R')^2=6$, a contradiction with Lemma \ref{ramif a}. On the other hand, since $K_S^2=8$, $S$ does not belong to case ($A4$), therefore $S$ belongs to case ($A3$), hence again Proposition \ref{Xiao list 2}, together with $\delta=0$, gives $(R')^2=10-2i$. Therefore Lemma \ref{ramif a} implies $i=5$, so Theorem \ref{Xiao list 1} shows that $S$ contains a rational pencil $|C|$ of hyperelliptic curves of genus $3$ with six double fibres. Now we have to prove that the pencil $|C|$ is base-point free. Suppose that $x \in S$ is a base point of $|C|$, and let $\sigma$ be the bicanonical involution of $S$; by Proposition \ref{relazioni}, $\sigma$ contains $12$ fixed points. If the image of $x$ in $S / \sigma$ is a smooth point of $S / \sigma$, then it gives rise to a base point of the ruling of $\widehat{W}$, and this is impossible. Then $x$ must be an isolated fixed point of $\sigma$, and therefore $\hS$ contains a $(-1)-$curve, pointwise fixed by $\hat{\sigma}$, which is a multisection of the pullback $|\widehat{C}|$ of $|C|$. But this is again impossible because, by Corollary \ref{fibre doppie}, the 12 $(-1)-$curves in $\widehat{S}$ come from the six $[5,5]-$points of $B$; then they are contained in fibres of $|\widehat{C}|$.\\  
It remains to show that the pencil $|C|$ is isotrivial. Consider the following base$-$change:
\begin{equation} 
\begin{CD}
X  @>\tilde{h}>> S\\
@VV{\tilde{\beta}}V  @VV{\beta}V\\
B @>h>> \mathbb{P}^ 1,\\
\end{CD}
\end{equation}
where $h$ is the double cover of $\mathbb{P}^1$ branched at the six
points corresponding to the six double fibres of $|C|$. It is clear
that $ \tilde{h}$ is an {\'e}tale double cover and that the fibres of
$\tilde{\beta}$ form on $X$ a genus $3$ pencil $|C'|$. Since
$K_X^2=16, \; g(B)=2$, we have $K_X^2=8(g(B)-1)(g(C')-1)$. Then the Arakelov theorem (see [Be82]) allows us to conclude that $\tilde{\beta}$ is an isotrivial fibration, and this implies that $\beta$ is an isotrivial fibration.  
\end{proof}

\subsection{The three families}

The next result shows that the existence of the isotrivial pencil $|C|$ gives very strong information about the geometry of $S$.

\begin{proposition} \label{struttura S}
Let $S$ be a minimal surface of general type with $p_g=q=1$, $K_S^2=8$ and bicanonical map of degree $2$. Then $S$ is a quotient of the form $(C \times F)/G$, where $C$, $F$ are smooth curves and $G$ is a finite group such that:
\begin{itemize}
\item $G$ acts faithfully on $C$, $F$, and freely on the product. Moreover, the action of $G$ on $C \times F$ is the diagonal one, namely $g(x,y)=(gx,gy)$.
\item $C$ is an hyperelliptic curve of genus $3$, and $C/G \cong E$, where $E$ is an elliptic curve isomorphic to the Albanese variety of $S$;
\item $F/G \cong \mathbb{P}^1$, and $|G|=2(g(F)-1)$.
\end{itemize}
\end{proposition}

\begin{proof}
Proposition \ref{isotrivial pencil} shows that $S$ contains an isotrivial pencil $|C|$ of hyperelliptic curves of genus $3$ such that the only singular fibres are twice a smooth curve of genus $2$. So following Serrano [Se90, Definition 1.2] we can say that $S$ is a \emph{quasi-bundle}, then [Se90, Theorem 2.1] shows that there exists a smooth curve $F$ and a finite group $G$ acting faithfully on $C$, $F$ and whose diagonal action is free on the product $C \times F$, in such a way that $S \cong (C \times F)/G$. The surface $S$ has two morphisms $S \longrightarrow F/G, \; S \longrightarrow C /G$, induced by the two projections of $C \times F$; notice that the fibres of the former morphism are the pencil $|C|$, hence we get $F/G \cong \mathbb{P}^1$. Moreover we have $q(S)=g(C/G)+g(F/G)$, so $C/G$ is an elliptic curve and this in turn means that the map $\alpha: S \longrightarrow C/G$ coincides with the Albanese fibration of $S$. Finally, recall that the invariants of $C \times F$ are 
\begin{displaymath}
\begin{split}
p_g(C \times F) & =g(C) \cdot g(F)=3g(F); \\
q(C \times F) & =g(C)+g(F)=3+g(F).
\end{split}
\end{displaymath}
Since $\pi: C \times F \longrightarrow S$ is an {\'e}tale Galois cover with Galois group $G$ and since $\chi(\mathcal{O}_S)=1$, we obtain $|G|= \chi(\mathcal{O}_{C \times F})=2(g(F)-1)$.  
\end{proof}

\begin{remark} \label{6 punti branched}
Since $G$ acts freely on the product $C \times F$, the singular elements of the two fibrations $\alpha$ and $\beta$ are multiple of smooth curves. In particular, since $\beta: S \longrightarrow \mathbb{P}^1$ contains six double fibres, it follows that the $G-$cover $f: F \longrightarrow \mathbb{P}^1$ is branched at six points $p_1, \ldots, p_6$, with branching order $2$ in each of them. 
\end{remark}
Now we want to understand which cases really occur for the pair $(g(F), G)$. The next theorem completely answers this question:

\begin{theorem} \label{gruppi}
Let $S$ be as in Proposition \ref{struttura S}.
We have exactly the following possibilities for $g(F)$ and $G$:
\begin{itemize}
\item[$I$.]  $g(F)=3$, $G=\mathbb{Z}_2 \times \mathbb{Z}_2$;
\item[$II$.]  $g(F)=4$, $G=S_3$;
\item[$III$.]  $g(F)=5$, $G=D_4$.
\end{itemize}
Moreover, each of these cases really occurs.
\end{theorem}

\begin{proof}
Here we show that $I, \;  II, \;  III$ are the only possibilities that can occur. The existence of the corresponding surfaces (that we shall call \emph{surfaces of type} $I, \; II, \; III$ respectively) will be shown in Section \ref{esempi}. \\ 
We consider again the two morphisms $\beta: S \longrightarrow F/G \cong \mathbb{P}^1$, $\alpha: S \longrightarrow C/G \cong E$. By the previous Remark \ref{6 punti branched}, the only singular fibres of $\alpha$ are multiple of smooth curves; moreover, a fibre of $\alpha$ appearing with multiplicity $n$ corresponds to a point $p \in E$ where the $G-$cover  $h:C \longrightarrow E$ is branched, and such that any point $q \in h^{-1}(p)$ has stabilizer group $G_q \cong \mathbb{Z}_n$. Now we apply the Zeuthen$-$Segre formula ([BPV84, Proposition 11.4 p.97]) to the fibration $\alpha: S \longrightarrow E$, or, equivalently, the Hurwitz formula to the covering $h: C\longrightarrow E$. If $\{ F_i=n_i \Phi_i \}_{i=1,\ldots,k}$ are the multiple fibres of $\alpha$, we obtain
\begin{displaymath}
4=c_2(S)=(2g(F)-2)\sum_{i=1}^k \Big( 1-\frac{1}{n_i} \Big),
\end{displaymath}
that is
\begin{equation} \label{formula}
2=(g(F)-1)\sum_{i=1}^k \Big(1-\frac{1}{n_i} \Big).
\end{equation}
Since $K_S^2=8$, Proposition \ref{caso standard} shows that we have $g(F) \geq 3$; moreover the sum in the right-hand side of (\ref{formula}) is clearly $\geq \frac{1}{2}$, hence $g(F) \leq 5$. Therefore we have only the following possibilities:
\begin{itemize}
\item[$I$.] $g(F)=3, \quad k=2, \quad n_1=n_2=2$;
\item[$II$.] $g(F)=4, \quad k=1, \quad n_1=3$;
\item[$III$.] $g(F)=5, \quad k=1, \quad n_1=2$.
\end{itemize}
We analyze separately each of these three cases.\\ \\
Case $I$.\\ 
The Albanese fibration of $S$ has genus $3$ and it contains two double fibres $F_1=2 \Phi_1$, $F_2=2 \Phi_2$, where $\Phi_i$ is a smooth curve of genus $2$. By Remark \ref{6 punti branched} the cover $f: F \longrightarrow \mathbb{P}^1$ is branched in $6$ points, and the stabilizer in each point of the ramification locus is isomorphic to $\mathbb{Z}_2$; on the other hand the cover $h: C \longrightarrow E$ is branched in two points, and the corresponding stabilizers are again isomorphic to $\mathbb{Z}_2$. By Proposition \ref{struttura S} the group $G$ has order $4$, hence we have only two possibilities: either $G=\mathbb{Z}_4$ or  $G=\mathbb{Z}_2 \times \mathbb{Z}_2$. But the former case is impossible, since $\mathbb{Z}_4$ contains just a subgroup of order $2$, hence its diagonal action on the product $C \times F$ cannot be free. It follows  $G=\mathbb{Z}_2 \times \mathbb{Z}_2$ in this case.\\ \\
Case $II$. \\ 
The Albanese fibration of $S$ has genus $4$ and it contains a triple
fibre $F_1=3 \Phi_1$, where $\Phi_1$ is a smooth curve of genus
$2$. The covering $h:C \longrightarrow E$ is branched at one point and
the two points in the corresponding fibre have stabilizer isomorphic
to $\mathbb{Z}_3$. By Proposition \ref{struttura S} the group $G$ has
order $6$, hence either $G=\mathbb{Z}_2 \times \mathbb{Z}_3$ or
$G=S_3$. But the former case must be excluded: indeed, by [Se90, Lemma
4.4] it follows that there does not exist any abelian cover of a
smooth algebraic curve branched in just one point. 
Therefore in this case we get $G=S_3$. \\ \\
Case $III$. \\
The Albanese fibration of $S$ has genus $5$ and it contains one double
fibre $F_1=2 \Phi_1$, where $\Phi_1$ is a smooth curve of genus
$3$. The covering $h:C \longrightarrow E$ is branched at one point,
and the four points in the corresponding fibre have stabilizer
isomorphic to $\mathbb{Z}_2$. By Proposition \ref{struttura S} it
follows that $|G|=8$, and the same argument as in the previous case allows us to conclude that $G$ is not abelian. So either $G =\mathbb{H}$ or $G=D_4$, where $\mathbb{H}=\{ \pm 1,\pm i, \pm j, \pm k \}$ is the group of quaternions, and $D_4$ is the dihedral group of order $8$. But the case $G=\mathbb{H}$ cannot occur, because $\mathbb{H}$ contains just one subgroup isomorphic to $\mathbb{Z}_2$, namely the one generated by $-1$, so its diagonal action cannot be free on the product. This shows that in this case $G=D_4$.  
\end{proof}
\begin{remark} \label{fissodivisoriale}
By Proposition \ref{relazioni}, the divisorial fixed locus $R'$ of the bicanonical involution $\sigma$ satisfies $K_SR'=8$. Now Lemma \ref{ramif a} and the analysis done in the proof of Theorem \ref{gruppi} allow us to describe $R'$ in the three cases. In particular we have:
\begin{itemize}
\item if $S$ is a surface of type $I$, then $R'=F_1+F_2$, where each $F_i$ is a smooth Albanese fibre; 
\item if $S$ is a surface of type $II$, then $R'=F + \Phi_1$, where $F$ is a smooth fibre and $\Phi_1$ is the (reduced) triple fibre;
\item if $S$ is a surface of type $III$, then $R'=F$, where $F$ is a smooth fibre.
\end{itemize}
\end{remark}

\begin{remark} \label{vale rat inv 2}
Looking at Remark \ref{vale rat inv} it follows that Proposition
\ref{struttura S} and Theorem \ref{gruppi} still hold if one does not require that the bicanonical map $\phi$ has degree $2$, but only that it factors through a rational involution.  
\end{remark}

\begin{definition}
We say that a surface $S$ is \emph{isogenous to a product} if there exists a finite {\'e}tale morphism $C \times F \longrightarrow S$, where $C, \; F$ are smooth curves. 
\end{definition} 
Notice that surfaces of type $I, \; II, \;III$ are surfaces isogenous to a product, since they are isomorphic to the quotient of a product of curves $C \times F$ by the free action of a finite group. Actually this is the general case, because Catanese in [Ca00] proved that if $S$ is a surface isogenous to a product, then there exist two smooth curves $C,\; F$ and a group $G$ acting freely on the product such that $S=(C \times F)/G$. Besides, he gave a topological characterization of surfaces isogenous to a product in terms of their fundamental group ([Ca00, Theorem 3.4]).

\section{Some technical tools} \label{tech tool}

\subsection{Automorphisms of hyperelliptic curves and Accola's Theorem} \label{hyp aut}

\subsubsection{Automorphisms of hyperelliptic curves}
It is well known that any smooth hyperelliptic curve $C$ of genus $g \geq 2$ is isomorphic to the smooth compactification of an affine curve $A$ of equation $y^2=p(x)$, where $p(x)$ is a polynomial without multiple roots of degree either $2g+1$ or $2g+2$, where the zeros of $p$ are ramifications points for the double cover $\pi: C \rightarrow \mathbb{P}^1$; the degree of $p$ is $2g+1$ if $\infty$ belongs to the ramification divisor of $\pi$ and $2g+2$ otherwise. Since the polynomial $p(x)$ is without multiple roots, the affine curve $A$ is smooth. The corresponding projective curve $\overline{A} \subset \mathbb{P}^2$ is singular at $\infty$, and $C$ will be the normalization of $\overline{A}$. Put $t= 1 / x$; therefore the point at infinity of $\overline{A}$ is defined by $t=0$. Let us consider the polynomial $k(t):=t^{2g+2}p(1/t)$. We have $k(t)=t^{2g+2}y^2(t)=(t^{g+1}y(t))^2$; hence setting $w=t^{g+1}y$ we obtain the affine curve $B$ defined by the equation $w^2=k(t)$. The hyperelliptic curve $C$ is therefore the quotient space
\begin{displaymath}
C=\frac{A \sqcup B}{\sim},
\end{displaymath}     
where the equivalence relation $\sim$ is given in the obvious way $(x,y) \sim (1/x, y(1/x)^{g+1})$; in this case, we shall say that the $C$ is defined by the equations
\begin{displaymath} 
 \left \{ \begin{array}{ll}
y^2=p(x) \\ 
w^2=k(t),
\end{array} \right.
\end{displaymath} 
or, more simply, that $C$ is defined by the equation $y^2=p(x)$.
Equivalently, $C$ is isomorphic to the curve $y^2=p(x_0,x_1)$ in the weighted projective space $\mathbb{P}(1,1,d)$, where $d:=g+1$ and $p(x_0,x_1)$ is a homogeneous polynomial of degree $2d$. Notice that the hyperelliptic involution $\tau$ is induced on $C$ by the automorphism $(x_0,x_1,y) \longrightarrow (x_0,x_1,-y)$, whilst the hyperelliptic double cover is the restriction to $C$ of the projection $\mathbb{P}(1,1,d) \longrightarrow \mathbb{P}^1$ given by $(x_0,x_1,y) \longrightarrow (x_0,x_1)$.\\ 
Now we consider the automorphism group of $C$, Aut($C$). It is always
different from the trivial group because $C$ possesses the hyperelliptic involution $\tau$. Since the $g^1_2$ on $C$ is unique, the automorphism $\tau$ is in the center of Aut($C$),  then for any $\phi \in$ Aut($C$) we have $\phi \circ \tau = \tau \circ \phi$. This implies that  the subgroup $ \langle \tau \rangle \cong \mathbb{Z}_2$ is normal in Aut($C$); moreover, if $G \subset$ Aut($C$) is any subgroup containing $\tau$, then $C /G \cong \mathbb{P}^1$.
Now let $\phi \in$ Aut($C$). The morphism $\pi \circ \phi$ has degree $2$; again by the uniqueness of the $g^1_2$ on $C$, there exists a projectivity $f$ which makes the following diagram commutative: 
\begin{equation} \label{diagram iperellittica}
\begin{CD}
C  @>\phi>> C\\
@VV{\pi}V  @VV{\pi}V\\
\mathbb{P}^1 @>f>> \mathbb{P}^1. \\
\end{CD}
\end{equation}
The projectivity $f$ acts on the polynomial $p$ multiplying it by a
non zero scalar $\lambda \in \mathbb{C}^{*}$; we will express this
fact by saying that $p$ is \emph{invariant} under $f$. This is
equivalent to saying that the zero locus $\Delta$ of $p(x_0,x_1)$ is
invariant under $f$. Vice versa, let us suppose that $f: \mathbb{P}^1
\rightarrow \mathbb{P}^1$ is a projectivity and let us try to lift it
to an automorphism $\phi$ of $C$. Of course, if $f$ admits a lifting
$\phi$, $\tau \circ \phi$ is also an allowed lifting; moreover, it is not difficult to see that there are no other possible liftings. Now suppose that $\phi$ does exist. Since the polynomial $p(x)$ is invariant under $f$, we must have $(p \circ f)(x)= \lambda p(x)$, $\lambda \in \mathbb{C}^*$. If the affine expression for $f$ is
\begin{displaymath}
f(x)= \frac{ax+b}{cx+d},
\end{displaymath}  
a straightforward computation (see [Pal03, Capitolo 4] ) shows that the two liftings of $f$ are
\begin{displaymath} 
\phi(x,y)=\Big ( \frac{ax+b}{cx+d}, \thinspace \frac{\sqrt{\lambda}y}{(cx+d)^{g+1}} \Big )  
\end{displaymath}
\begin{displaymath} 
\tau \circ \phi(x,y)= \Big( \frac{ax+b}{cx+d}, \thinspace - \frac{\sqrt{\lambda}y}{(cx+d)^{g+1}} \Big).
\end{displaymath}
Now suppose that a finite group $G$ acts on $C$. Then the above
argument shows that the action of $G$ descends to an action of a group
$\widehat{G}$ on $\mathbb{P}^1=C/ \tau$ which preserves the branch
locus $\Delta$ of the double cover $C \longrightarrow
\mathbb{P}^1$. If we identify $\mathbb{P}^1$ with the unitary sphere
in $\mathbb{R}^3$, we obtain an inclusion $SO(3)
\hookrightarrow$Aut$(\mathbb{P}^1)$; moreover, by a conjugation
argument it is possible to see that all the finite subgroups of
Aut$(\mathbb{P}^1)$ can be obtained from a finite subgroup of $SO(3)$ in this way.  On the other hand, the finite subgroups of $SO(3)$ are just the following:
\begin{itemize}
\item the cyclic groups $\mathbb{Z}_n$;
\item the dihedral groups $D_n$ of order $2n$, $n \geq 2$ (where $D_2=\mathbb{Z}_2 \times \mathbb{Z}_2$);
\item the group of symmetries of the tetrahedron, which is isomorphic to $A_4$;
\item the group of symmetries of the cube (and of the octahedron) which is isomorphic to $S_4$;
\item the group of symmetries of the dodecahedron (and of the icosahedron) which is isomorphic to $A_5$.
\end{itemize}
Notice moreover that, denoting by $G_0$ the subgroup of Aut$(C)$ generated by $ \langle \tau \rangle$ and $G$, we have a central extension
\begin{equation} \label{estensione centrale}
0 \longrightarrow  \langle \tau \rangle \longrightarrow G_0 \longrightarrow \widehat{G} \longrightarrow 0.
\end{equation}
If $\tau \in G$, then $G=G_0$, whereas if $\tau \notin G$ (and this
will always be our case) then $G$ is mapped homomorphically onto
$\widehat{G}$. Let $q: SL(2,\mathbb{C}) \longrightarrow
PGL(1,\mathbb{C})$ be the quotient map and let $H:=q^{-1} \widehat{G}$. We have a surjective map $H \longrightarrow \widehat{G}$ whose kernel is the subgroup generated by $-$Id. The group $H$ acts on the polynomials of degree $2g+2$ and $p(x_0,x_1)$ is an eigenvector for this action. So we have a homomorphism $\lambda: H \longrightarrow \mathbb{C}^*$ defined by $(h^{-1})^*p(x_0,x_1)= \lambda(h)p(x_0,x_1)$. But the degree of $p(x_0,x_1)$ is even, so $-$Id $\in$ Ker($\lambda$) and we have actually defined a character $\lambda : \widehat{G} \longrightarrow \mathbb{C}^*$. In all our examples the character $\lambda$ will be trivial, so we will apply the following result; see [Par03, p.101].  

\begin{proposition} \label{carattere banale}
Suppose that $\lambda$ is the trivial character, and that $g$ is odd. Then the exact sequence $($\ref{estensione centrale}$)$ splits. If moreover  $\tau \notin G$, then $G_0 \cong G \times \mathbb{Z}_2$.
\end{proposition}
\begin{proof}
Let us consider the homomorphism $H \longrightarrow G_0$ given in the following way: we associate to $h \in H$ the automorphism of $C$ given by $(x_0,x_1,y) \longrightarrow (h(x_0,x_1),y)$. As $g$ is odd, $d=g+1$ is even and this in turn implies that $-$Id belongs to the kernel of $H \longrightarrow G_0$; so we have actually constructed a homomorphism $\widehat{G} \longrightarrow G_0$ that splits the central extension (\ref{estensione centrale}).  Then $G_0$ is the semidirect product of $\langle \tau \rangle$ and $\widehat{G}$, and since the only automorphism of $\mathbb{Z}_2$ is the identity, $G_0$ is actually a direct product, i.e. $G_0 \cong \widehat{G} \times \mathbb{Z}_2$. If moreover $\tau \notin G$, then $\widehat{G} \cong G$, so $G_0 \cong G \times \mathbb{Z}_2$ and we are done.
\end{proof}

\subsubsection{Accola's theorem}

\begin{definition} \label{gruppo con partiz}
A finite group $G_0$ is said to admit a \emph{partition} if there is a collection of subgroups $$G_1, G_2, \ldots ,G_s \; \; \; (s \geq 3)$$ such that:
\begin{itemize}
\item $G_0= \bigcup_{i=1}^s G_i$;
\item $G_i \cap G_j =id  \;$ for $\;0 <i<j$.
\end{itemize}
\end{definition}
Automorphism groups of Riemann surfaces admitting a partition are interesting because there is the following Theorem \ref{accola}.

\begin{theorem}[$\mathbf{Accola}$] \label{accola}
Let $C$ be a smooth curve of genus $g$, and let $G_0$ be a subgroup of Aut$(C)$ admitting a partition $G_1, \ldots, G_s$.
Let $n_i$ be the order of $G_i$ and let $g_i$ be the genus of the curve $C/G_i$. Then
\begin{equation} \label{accolaformula}
(s-1)g+n_0g_0=\sum_{i=1}^sn_ig_i.
\end{equation}
\end{theorem}
\begin{proof}
See [Ac94, Theorem 5.9].
\end{proof}
Formula (\ref{accolaformula}) is of interest because the ramifications
do not enter into it. Now we want to look at the theorem of Accola in
two particular cases: when $G_0=(\mathbb{Z}_2)^n$ and when $G_0=D_n$ (the dihedral group of order $2n$). \\ \\
$\bullet \;  \; G_0=(\mathbb{Z}_2)^n$. \\  
In this case $G_0$ admits a partition of type $G_0= \bigcup_{i=1}^{2^n-1}G_i$, where $G_i \cong \mathbb{Z}_2$; then formula (\ref{accolaformula}) becomes
\begin{equation} \label{accolaformula1}
(2^{n-1}-1)g+2^{n-1}g_0= \sum_{i=1}^{2^n-1}g_i.
\end{equation}
$\bullet \; \; G_0=D_n$. \\
In this case $G_0$ is generated by a rotation $r$, which has order
$n$, and by a reflection $s$, which has order $2$. Let $G_r \cong
\mathbb{Z}_n$ be the subgroup of $G_0$ generated by $r$ and let $G_i \cong \mathbb{Z}_2$ be the subgroup generated by the element $r^i s$, $i=1, \ldots, n$. Then $G_0$ admits the partition
\begin{displaymath}
G_0=G_r \cup \big( \bigcup_{i=1}^n G_i \big).
\end{displaymath}
Set $C_r:=C/G_r$, $C_i:= C/G_i$, $g_r:=g(C_r)$ and $g_i:=g(C_i)$. Theorem \ref{accola} gives
\begin{equation} \label{accolaformula2}
ng+2n g_0=ng_r+ 2 \sum_{i=1}^ng_i.
\end{equation}
 If $n$ is odd, then all the subgroups $G_i$ are conjugated, and so $g_i=g_1$, $i=1, \ldots, n$; in this case formula (\ref{accolaformula2}) becomes
\begin{equation} \label{accolaformula3a}
g+2 g_0=g_r+2g_1.
\end{equation}
If $n$ is even, then $G_i$ is conjugate to $G_j$ if and only if $i \equiv j$ (mod $2$); in this case formula (\ref{accolaformula2}) becomes
\begin{equation} \label{accolaformula3b}
g + 2 g_0= g_r +g_1+g_2. 
\end{equation} \\
We say that a curve $C$ is \emph{bielliptic} if it is a double cover of an elliptic curve. Accola's Theorem allows us to prove the following interesting result (see also [BaDC99, Theorem 2.1]):

\begin{theorem} \label{iperellbiell}
Let $C$ be a curve of genus $3$. Then the following are equivalent:
\begin{itemize}
\item[($1$)] $C$ is a double cover of a curve of genus $2$.
\item[($2$)] $C$ is both hyperelliptic and bielliptic.
\end{itemize}
\end{theorem}

\begin{proof}
($1$) $\Longrightarrow$ ($2$). \\
Assume that there exists a double cover $f: C \longrightarrow D$, where $D$ is a smooth curve of genus $2$. Notice that $f$ is an  \' etale morphism (by Hurwitz formula). Let $ \sigma : C \longrightarrow C$ be the involution of $C$ such that $D=C / \sigma$, and let $T: D \longrightarrow D$ be the hyperelliptic involution of $D$. We claim that $T$ lifts to an involution $\tau$ of $C$. Indeed, observe that if $x, \thinspace y$ are points of $D$, then $x+T(x)$ is linearly equivalent to $y+T (y)$, hence the involution $T$ induces the automorphism $p \rightarrow -p$ on the Jacobian $J(D)$ of $D$, and this in turn implies that if $\eta \in$ Pic($D$) is an element of $2-$torsion, then $T_* \eta \cong \eta$. Since the double covering $f$ is {\'e}tale, it is defined by a $2-$torsion element $\eta$, so the involution $T$ lifts to an automorphism $\tau \in $ Aut($C$) which makes the following diagram commutative:

\begin{equation} \label{sollevamento T}
\begin{CD}
C  @>\tau>> C\\
@VV fV     @VVf V\\
D @>T>> D. \\
\end{CD}
\end{equation}
The order of $\tau$ is either $2$ or $4$. If $\tau$ has order $4$, then $ \langle \tau \rangle \cong \mathbb{Z}_4$, and diagram (\ref{sollevamento T}) shows $\tau^2=\sigma$. Hence $\tau^2$ does not have fixed points on $C$, and so the same holds for $\tau$ and $\tau^3$ ( notice that $(\tau^3)^2=\tau^2$ ). Therefore the quadruple covering $C \longrightarrow C/ \langle \tau \rangle$ is {\'e}tale; if $g$ is the genus of $C /\langle \tau \rangle$, Hurwitz formula gives $4=4(2g-2)$, thus $2g-2=1$, a contradiction. So $ \tau$ is an involution, and this proves our claim. \\
Now, from diagram (\ref{sollevamento T}) it follows easily that $\tau \sigma = \sigma \tau$, and this means that the group $G_0= \{id, \sigma, \tau, \sigma \tau \}$ is isomorphic to $\mathbb{Z}_2 \times \mathbb{Z}_2$; moreover, since $\tau$ is a lifting of the hyperelliptic involution $T$ of $D$, we get $C/G_0 \cong D /T \cong \mathbb{P}^1$. Hence if we set $g_{\tau}:=g(C/\tau)$, $g_{\sigma \tau}:= g(C/\sigma \tau)$, formula (\ref{accolaformula1}) becomes
\begin{displaymath}
3=2+g_\tau +g _{\sigma \tau}.
\end{displaymath}
So we have $\{ g_\tau, g_{\sigma\tau} \}= \{0,1 \}$, and this shows that $C$ is both  hyperelliptic and bielliptic.\\  \\
($2$) $\Longrightarrow$ ($1$). \\
Assume that $C$ is hyperelliptic and bielliptic, and let $\tau$ be the hyperelliptic involution, and $\sigma$ be the involution induced by the bielliptic structure. By definition $\tau ^2= \sigma ^2=1$, and since $ \tau $ is a central element of Aut($C$), we have $\sigma \tau = \tau \sigma$. Hence $ \{ id, \tau, \sigma, \tau\sigma \}$ is a subgroup of Aut($C$) isomorphic to $\mathbb{Z}_2 \times \mathbb{Z}_2$. Applying formula (\ref{accolaformula1}) we obtain
\begin{displaymath}
3=1+g(C/ \tau \sigma),
\end{displaymath}
that is, $g(C /\tau \sigma)=2$, and we are done.
\end{proof}

\begin{remark} \label{biell no 4}
Theorem \ref{iperellbiell} cannot be generalized to curves of higher genus. Indeed, a bielliptic curve $C$ of genus $g(C)\geq 4$ is not hyperelliptic. This is a consequence of the so-called Castelnuovo-Severi inequality; see for example $[Ac94]$, $[$ACGH85, Exercise C-1 p.366$]$ and $[$Xi87, Lemma 7 p.465$]$.
\end{remark} 

\subsection{Orbifold fundamental group and Hurwitz monodromy}

\subsubsection{The orbifold fundamental group}
Now we introduce the orbifold fundamental group, which is a powerful
tool used to study the Galois covers of a complex manifold. Since we are mostly interested in Galois covers of algebraic curves, we will present the theory only in this case, and we refer the reader to [Ca00] for a more complete exposition. Let $f: X \longrightarrow Y:=X /G$ be the quotient of a smooth algebraic curve $X$ with respect to a finite group of automorphisms. Then $Y$ is again a smooth algebraic curve, and we denote by $B= \{ p_1, \ldots, p_r \}$ the set of points in $Y$ where $f$ is branched. The ramification locus of $f$ is defined as the locus of points $q \in X$ such that $f(q)=p_i$ for some $i$. If $m_i$ is the branching order of $f$ in $p_i$, then the fibre $f^*(p_i)$ consists of $d=\frac{|G|}{m_i}$ points $q_1, \ldots, q_d$, and the stabilizer of any $q_j$ is a cyclic subgroup of $G$ isomorphic to $\mathbb{Z}_{m_i}$. Besides, the stabilizer subgroups of two points in the same fibre of $f$ are conjugated in $G$. Let $x_i \in \pi_1(Y -B)$ be a simple geometric loop around $p_i$; we have an exact sequence
\begin{displaymath} 
1 \longrightarrow \pi_1(X- f^{-1}(B)) \longrightarrow \pi_1(Y -B) \longrightarrow G \longrightarrow 1,
\end{displaymath} 
and each $x_i$ maps to an element $g_i \in G$ of order $m_i$; let $\mathbf{m}=(m_1, \ldots, m_r)$.

\begin{definition}
The \emph{orbifold fundamental group} $\pi_1^{orb}(Y -B, \; \mathbf{m}, \; y_0)$ is the quotient of $\pi_1(Y -B, \; y_0)$ by the minimal normal subgroup generated by the elements $(x_i)^{m_i}$.
\end{definition}
Where there is no possibility of confusion, we will omit the base point $y_0$ and we will simply write $\pi_1^{orb}(Y -B, \; \mathbf{m})$. \\    
If the genus of $Y$ is $g$, we can think of $Y$ as the quotient space obtained by identifying in a suitable way the edges of a regular $4g-$gon. If $x_1, \ldots, x_r$ are simple loops around the points $p_i$ in $Y -\{ p_1, \ldots, p_r \}$, then $\pi_1^{orb}(Y - \{ p_1, \ldots, p_r \}, \; \mathbf{m})$ is the Fuchsian group admitting the presentation
\begin{equation} \label{fuchsian} 
 \langle x_1, \ldots, x_r; \; a_1, \ldots, a_g, \; b_1, \ldots ,b_g \; | \; x_1^{m_1}=x_2^{m_2}= \ldots =x_r^{m_r}=x_1x_2 \cdots x_r \cdot \prod_{i=1}^g [a_i,b_i]=1 \rangle.
\end{equation}

\begin{definition}
An \emph{admissible epimorphism} $\mu: \pi_1^{orb}( Y - B, \; \mathbf{m}) \longrightarrow G$ is a surjective homomorphism of groups such that $g_i:=\mu(x_i)$ is an element of $G$ of order $m_i$, for any $i=1, \ldots, r$.    
\end{definition}
What follows is the generalization to branched coverings of a result which is well known in the unramified case.
\begin{theorem} \label{orbifold fond}
Let $Y$ be an algebraic curve, and let $G$ be a finite group. Let $B=\{p_1, \ldots, p_r \}$ be a set of $r$ distinct points on $Y$, and let $\mathbf{m}:=(m_1, \ldots, m_r)$ be an $r-$ple of positive integers. Then there exists a ramified Galois cover $f:X \longrightarrow Y$ branched at $B$ in such a way that:
\begin{itemize}
\item the Galois group of $f$ is $G$;
\item the branching order of $f$ at $p_i$ is $m_i$
\end{itemize}
if and only if there exists a admissible epimorphism
\begin{displaymath} 
\mu: \pi_1^{orb} (Y-B, \; \mathbf{m}) \longrightarrow G.
\end{displaymath}
Moreover, in this case one has a short exact sequence of groups:
\begin{displaymath} 
1 \longrightarrow \pi_1(X) \longrightarrow \pi_1^{orb} (Y-B, \; \mathbf{m}) \xrightarrow{\; \; \mu \;} G \longrightarrow 1,
\end{displaymath}
and this shows that the isomorphism class of $X$ is determined not by
the specific epimorphism $\mu$, but rather by its kernel.
\end{theorem}

\begin{remark}
For any $g \in G$, we denote by $\textrm{Fix}(g):=\{x \in X \;| \; gx=x \}$ the set of fixed points of $g$. $\textrm{Fix}(g) \neq \emptyset$ if and only if the cyclic subgroup $ \langle g \rangle$ generated by $g$ is contained in the stabilizer of some point $q$ in the ramification locus of $f$. It is easy to see that these stabilizers are exactly the subgroups $\langle g_i \rangle$ and all their conjugates in $G$.
\end{remark}

\subsubsection{Hurwitz monodromy}
Now let us consider the two spaces
\begin{displaymath}
\begin{split}
& F_rY:= \{ (y_1, \ldots, y_r) \; | \;  y_i=y_j \; \; \textrm{iff} \; \; i=j \} \subset Y^r; \\
&B_rY:= \{ y_1+ \cdots +y_r \; | \;  y_i=y_j \; \; \textrm{iff} \; \; i=j \} \subset \textrm{Sym}^r \; Y.
\end{split}
\end{displaymath}

\begin{definition} 
$\pi_1(F_rY)$ is called the \emph{pure braid group} of $Y$, whereas $\pi_1(B_rY)$ is called the \emph{full braid group} of $Y$, or simply the \emph{braid group}.
\end{definition}
Notice that there is a Galois {\'e}tale covering $F_rY \longrightarrow B_rY$ with Galois group $S_r$, which gives rise to a short exact sequence:
\begin{displaymath}
1 \longrightarrow \pi_1(F_rY) \longrightarrow \pi_1(B_rY) \longrightarrow S_r \longrightarrow 1.
\end{displaymath}
Let us keep the curve $Y$ fixed; the set of the isomorphism classes of
Galois $G-$covers $f: X \longrightarrow Y$, branched in $r$ distinct
points with branching orders $m_1, \ldots, m_r$, will be called a \emph{Hurwitz space} and will be denoted by $\mathcal{H}(Y)^{G, \mathbf{m}}$. Let us consider the map
\begin{displaymath} 
\pi: \mathcal{H}(Y)^{G, \mathbf{m}} \longrightarrow B_rY
\end{displaymath}
given by associating to the $G-$cover $f: X \longrightarrow Y$ its
branch locus. Theorem \ref{orbifold fond} implies that
\begin{displaymath}
\pi^{-1}(B)= \{ \textrm{kernels of admissible epimorphisms } \mu: \pi_1^{orb}(Y -B, \; \mathbf{m}) \longrightarrow G \}.
\end{displaymath}
It is clear that the cardinality of $\pi^{-1}(B)$ does not depend on $B$, because the number of the admissible epimorphisms $\mu$ does not depend on the choice of the points $p_1, \ldots, p_r$, if such points remain distinct. It is actually possible give a topology to $\mathcal{H}(Y)^{G, \mathbf{m}}$ in such a way that $\pi$ becomes an unramified cover; see[BaCa97, p.421].
Choose a base point $B_0$ in $B_rY$; the monodromy
\begin{displaymath}
T_{\pi}: \pi_1(B_rY, \; B_0) \longrightarrow \{ \textrm{permutations of }\pi^{-1}(B_0) \}
\end{displaymath}
of the cover $\pi$ is called the \emph{Hurwitz monodromy}. Clearly the
topological space $\mathcal{H}(Y)^{G, \mathbf{m}}$ is connected if and
only if the Hurwitz monodromy is transitive. Consider the universal
family $u: \mathcal{U} \longrightarrow B_rY$ over $B_rY$, where
$\mathcal{U}=(B_rY \times Y)- \{(B,x)\; | \; x  \in B \}$ (roughly
speaking, the fibre of $u$ over $B \in B_rY$ is the punctured manifold
$Y-B$). Let $B_0 \in B_rY$, let $\gamma: I=[0,1] \longrightarrow B_rY$
be a  loop with $\gamma(0)=\gamma(1)=B_0$ and consider the pullback
$\gamma^*\mathcal{U}$ of $\mathcal{U}$ to $I$. Since $I$ is
contractible, the fibre bundle $\gamma^* \mathcal{U}$ is trivial; then
there is a bundle diffeomorphism $\psi: \gamma^* \mathcal{U}
\longrightarrow I \times (Y -B_0)$. $\psi_0$ is the identity map,
whereas $\psi_1: Y-B_0 \longrightarrow Y - B_0$ is the image of
$\gamma$ via the monodromy map $T_u: \pi_1(B_rY, \; B_0)
\longrightarrow \textrm{Diff}\;(u^{-1}(B_0))$ of $u$. For sake of
simplicity let us suppose that there is a subset $V \subset B_rY$ such
that there exists a section $s: V \longrightarrow \mathcal{U}$ of $u$;
this will always be true in our applications. Then we can choose the base point $y_0$ in such a way that $\psi_1(y_0)=y_0$. Let us consider the automorphism $\pi_1(Y-B_0, \; y_0) \longrightarrow \pi_1(Y-B_0, \; y_0)$ associated to $\psi_1$, and let $(\psi_1)_*$ be the induced automorphism of $\pi_1^{orb}(Y-B_0, \; \mathbf{m})$. What follows is the analogous of [BaCa97, Proposition 1.12] in the case of $G-$covers:

\begin{proposition} \label{mon p}
The monodromy of the cover $\pi: \mathcal{H}(Y)^{G, \mathbf{m}} \longrightarrow B_rY$ is described in the following way. $\pi_1(B_rY)$ acts as $\textrm{ker} \; (\mu) \longrightarrow \textrm{ker} \; (\mu \circ (\psi_i)_*^{-1})$ on the set
\begin{displaymath}
\pi^{-1}(B)= \{ \textrm{kernels of admissible epimorphisms } \mu: \pi_1^{orb}(Y -B, \; \mathbf{m}) \longrightarrow G \}.
\end{displaymath}
\end{proposition}
In what follows we are interested in applying this theory to two particular cases: when $Y \cong  \mathbb{P}^1$ and when $Y \cong E$, where $E$ is an elliptic curve. \\ \\ 
$\mathbf{Y = P^1.}$\, We denote $\pi_1^{orb}(\mathbb{P}^1-\{p_1, \ldots, p_r \}, \; \mathbf{m})$ by $\Gamma(m_1, \ldots,m_r)$. Let $x_i$ be a simple loop in $\mathbb{P}^1- \{p_1, \ldots, p_r \}$ around $p_i$. Then a presentation of  $\Gamma(m_1, \ldots, m_r)$ is
\begin{displaymath}   
\langle \; x_1, \ldots ,x_r \; | \; x_1^{m_1}= \ldots =x_r^{m_r}= x_1x_2 \cdots x_r=1 \; \rangle. 
\end{displaymath}
Let us consider a Galois cover $f: X \longrightarrow \mathbb{P}^1$ branched at $p_1, \ldots, p_r$ with branching orders $m_1, \ldots, m_r$; then by Theorem \ref{orbifold fond} we have an exact sequence
\begin{equation} \label{exact 3}
1 \longrightarrow \pi_1(X) \longrightarrow \Gamma(m_1, \ldots, m_r) \longrightarrow G \longrightarrow 1,
\end{equation}
 and the image of $x_i$ in $G$ is an element $g_i$ whose order is exactly $m_i$. Then the Galois cover $X$ exists if and only if it is possible to find elements $g_1, \ldots, g_r \in G$ of order $m_1, \ldots, m_r$ such that the $g_i$'s generate $G$ and $g_1g_2 \cdots g_r=1$. \\ 
The group $\pi_1(B_r\mathbb{P}^1)$ admits a presentation with generators $\sigma_1, \ldots, \sigma_{r-1}$ and defining relations (see [Bir74, p.34])
\begin{equation} \label{braid sphere}
\begin{split}
&\sigma_i\sigma_j=\sigma_j \sigma_i \quad \textrm{if }|i-j|\geq2 \\
& \sigma_i\sigma_{i+1}\sigma_i=\sigma_{i+1} \sigma_i \sigma_{i+1} \\
& \sigma_1\cdots \sigma_{r-2} \; \sigma_{r-1}^2 \sigma_{r-2} \cdots \sigma_1=1.
\end{split}
\end{equation}
It is not difficult to write down explicitly the monodromy action of $\pi_1(B_r\mathbb{P}^1)$. The automorphism $(\psi_1)_*$ of $\Gamma(m_1, \; \ldots ,m_r)$ corresponding to the generator $\sigma_i$ is simply the following: 
\begin{equation*} 
\begin{split}
& x_1 \longrightarrow x_1 \\
& \ldots \\
& x_{i-1} \longrightarrow x_{i-1} \\
& x_i  \longrightarrow x_ix_{i+1}x_i^{-1} \\
& x_{i+1}  \longrightarrow x_i \\
& x_{i+2}   \longrightarrow x_{i+2} \\
& \ldots \\ 
& x_r \longrightarrow x_r
\end{split}
\end{equation*}
(see [Bir74, p.25] and [Hur891]). We will call it the \emph{Hurwitz move} among $x_i$ and $x_{i+1}$. The \emph{Hurwitz equivalence} will be the equivalence relation on the set 
\begin{equation*} 
\{ \mu: \Gamma(m_1, \ldots, m_r) \longrightarrow G \; | \; \mu \textrm{ is an admissible epimorphism} \}
\end{equation*}
generated by:
\begin{itemize}
\item composition with automorphisms of $G$;
\item composition with Hurwitz moves.
\end{itemize}
Notice  that if the group $G$ is abelian, then a Hurwitz move is
simply a permutation of the generators $x_i$ of $\Gamma(m_1, \ldots,
m_r)$. An immediate consequence of Proposition \ref{mon p} is the following
\begin{proposition} \label{componenti di Hurwitz}
Let $\mu_1: \Gamma(m_1, \ldots ,m_r) \longrightarrow G$, $\mu_2:  \Gamma(m_1, \ldots ,m_r) \longrightarrow G$ be two admissible epimorphisms. Then the induced $G-$coverings $f_1: X_1 \longrightarrow \mathbb{P}^1$, $f_2: X_2 \longrightarrow \mathbb{P}^1$ belong to the same connected component of the Hurwitz scheme $\mathcal{H}(\mathbb{P}^1)^{G, \mathbf{m}}$ if and only if $\mu_1$ and $\mu_2$ are Hurwitz equivalent.
\end{proposition}
 \noindent $\mathbf{Y=E,\; \;E \;\;an \;\; elliptic\;\; curve.}$ \, In this case we denote $\pi_1^{orb}(E- \{p_1, \ldots, p_r \}, \; \mathbf{m})$ by $\Delta(m_1, \ldots, m_r)$. If $x_1, \ldots,x_r$ are simple loops in $\pi_1(E- \{ p_1, \ldots, p_r \})$ around $p_1, \ldots, p_r$, then a presentation of $\Delta(m_1, \ldots, m_r)$ is the following:
\begin{displaymath}
\langle \; x_1, \ldots, x_r; \; a, \; b \; | \; x_1^{m_1}= \cdots = x_r^{m_r}=x_1x_2 \cdots x_r [a,b] =1 \; \rangle.
\end{displaymath}  
Suppose that we have a Galois covering $f: X \longrightarrow E$ with Galois group $G$, branched at $p_1, \ldots ,p_r$ with branching orders $m_1, \ldots ,m_r$. Then we can write down the exact sequence
\begin{displaymath} 
1 \longrightarrow \pi_1(X) \longrightarrow \Delta(m_1, \ldots, m_r) \longrightarrow G \longrightarrow 1,
\end{displaymath}
and the image of  $x_i$ is an element $g_i \in G$ of order $m_i$. In particular, if $r=1$ we have
\begin{displaymath} 
\Delta(m)=\langle \; a,b \; | \; [a,b]^m=1 \; \rangle,
\end{displaymath}
and an admissible epimorphism $\delta: \Delta(m) \longrightarrow G$ is given by sending $a, \; b$ into two elements $g_1, \; g_2 \in G$ such that $g_1, \; g_2$ generate $G$ and the order of $[g_1, \; g_2]$ is $m$.

\section{The construction of the examples}  \label{esempi}
In this section we give examples of surfaces $S$ of type $I, \; II, \; III$, showing that they actually exist. This will be done by  constructing the curves $F \;, C$ and the action of the group $G$. We prove that the curve $C$ is both hyperelliptic and bielliptic, by using Accola's Theorem \ref{accola}. Moreover, in each case we write down an affine model for $C$ and we describe the action of $G$ explicitly, following the ideas contained in [Par03] and [Pal03]. It will turn out that both $C$ and $F$ belong to connected families. For $C$ this is a consequence of its explicit equation (in fact the family which parameterizes $C$ is also irreducible), whereas for $F$ this follows from Proposition \ref{componenti di Hurwitz}. 
  
\subsection{Surfaces of type $I$} \label{surfaces of type I}
In this case $G=\mathbb{Z}_2 \times \mathbb{Z}_2$. \\
We first construct the curve $C$. Let $E$ be an elliptic curve and let
$D=p+q$ be an effective divisor of degree $2$ over $E$, with $p$ and $q$ distinct points.
Let $D_1=D=p+q$, $D_2=D_3=0$. Since it is possible to find three elements $L_1, L_2, L_3$ in $\textrm{Pic}(E)$, $L_i \neq 0$ satisfying the relations $2L_i=D_j+D_k$, $L_i+L_j=D_k+L_k$, by [Ca84, Proposition 2.3] it follows that there exists a smooth $G-$cover $h: C \longrightarrow E$ branched over $D$. Then Hurwitz formula gives $g(C)=3$. Moreover $h$ factors through the double covers $h_{ij}: C_{ij} \longrightarrow E$ branched over $D_i+D_j$ and determined by the square root $L_k$. Hence we have the following commutative diagram:
\begin{equation}
\label{rombo}
\xymatrix@!0{
& & C \ar[lldd]_{g_{12}} \ar[dd]^{g_{13}} \ar[rrdd]^{g_{23}}  & & \\
& & & & & \\
C_{12} \ar[rrdd]_{h_{12}} & & \ C_{13} \ar[dd]^{h_{13}} & & C_{23} \ar[lldd]^{h_{23}} \\
& & & & &  \\
& & E & & \\
}
\end{equation}
where $g(C_{12})=2, \ g(C_{13})=2, \ g(C_{23})=1$.
Thus $C$ is a double cover of a smooth curve of genus $2$, so it is both hyperelliptic and bielliptic by Theorem \ref{iperellbiell}. \\ 
Now we will give an explicit equation for the curve $C$. Diagram (\ref{rombo}) shows that $\tau \notin G$, because there is no involution in $G$ with rational quotient, hence using the notations of Subsection \ref{hyp aut} we have $\widehat{G}=G$, and the action of $G$ on $C$ descends to an action of $G$ on $\mathbb{P}^1$ given in the following way: if $e_1, \; e_2$ are standard generators for $G$ and $e_3=e_1+e_2$, then
\begin{equation} \label{G1}
\begin{split}
e_1(x_0, x_1) & =(x_0,-x_1) \\
e_2(x_0, x_1) & =(x_1, x_0) \\
e_3(x_0,x_1)&=(x_1,-x_0).
\end{split}
\end{equation}
Thus we have to look for a polynomial $p(x)$ of degree $8$ with no multiple roots such that its zero locus $\Delta$ is invariant under  (\ref{G1}); if we write $p(x)= \prod_{i=1}^8 (x-\alpha_i)$, $\Delta$ being invariant under the action of $G$ means that the set $\{ \alpha_1, \ldots, \alpha_8 \}$ of zeroes of $p$ is of type $\{ \alpha, \frac{1}{\alpha}, \beta, \frac{1}{\beta}, -\alpha, -\frac{1}{\alpha}, -\beta, -\frac{1}{\beta} \}$. Then if $a:=\alpha^2, \; b:=\beta^2$, the general form of $p(x)$ is $p(x)=(x^2-a)(x^2-\frac{1}{a})(x^2-b)(x^2-\frac{1}{b})$. Therefore the equation of the curve $C$ is 
\begin{equation} \label{C prima}
y^2=(x^2-a)(x^2-b) \Big( x^2-\frac{1}{a} \Big) \Big( x^2-\frac{1}{b} \Big), 
\end{equation}   
where $a, \;b$ are such that the zeroes of $p(x)$ are distinct. Of course (\ref{C prima}) shows that $C$ belongs to an irreducible, $2-$dimensional family. \\
Lift the action of $G$ on $C$ as follows:
\begin{displaymath}
\begin{split}
e_1(x,y)&=(-x,y) \\
e_2(x,y) & =(1/x,-y/x^4) \\
e_3(x,y) &=(-1/x,-y/x^4). \\
\end{split} 
\end{displaymath}
One can immediately check that $e_2$ and $e_3$ have no fixed points, whereas $e_1$ has the four fixed points $(0,1), \thinspace  (0,-1), \thinspace ( \infty, 1), \thinspace (\infty, -1)$. Notice that the only other possible lifting does not work for our purposes, because it gives rise to an action which does not have the desired stabilizers. \\
Now we have to construct $F$. It will be defined by an admissible epimorphism
\begin{equation}
\mu: \Gamma(2,2,2,2,2,2) \longrightarrow G,  
\end{equation}
and the fact that $G$ acts freely on $C \times F$ means that $\mu$ must satisfy the following condition:
\begin{itemize}
\item[($*$)] the images of the generators $x_1, \ldots, x_6$ of $\Gamma(2,2,2,2,2,2)$ are different from $e_1$. 
\end{itemize}
An example of such a $\mu$ is, for instance
\begin{equation}
\begin{split}
x_1, \; x_2, \; x_3, \; x_4  & \longrightarrow e_2 \\
x_5, \; x_6        & \longrightarrow e_3,
\end{split}
\end{equation}
and it is clear that this is the unique possibility modulo automorphisms of $G$ and Hurwitz moves (namely, permutations of the $x_i$'s being $G$ abelian). Therefore Proposition \ref{componenti di Hurwitz} says that in this case the curve $F$ belongs to a connected family; this family has dimension $3$, because we are choosing six points on $\mathbb{P}^1$. \\
Notice that the $G-$cover $f: F \longrightarrow \mathbb{P}^1$ factors in the following way:
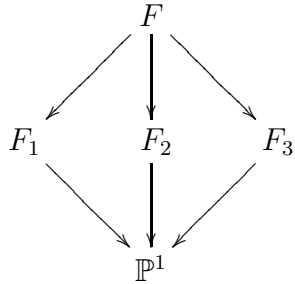
\[
\xymatrix@!0{
& & F \ar[lldd]  \ar[dd]  \ar[rrdd]  & & \\
& & & & & \\
F_1 \ar[rrdd]  & & \ F_2 \ar[dd]  & & F_3 \ar[lldd] \\
& & & & &  \\
& & \mathbb{P}^1 & & \\
}
\]
where all the arrows are double coverings, in fact $F_i=F/\langle e_i \rangle$. Observe that $g(F_1)=2, \ g(F_2)=0, \ g(F_3)=1$, in particular $F$ is \emph{hyperelliptic}. By construction the diagonal action of $G$ on $C \times F$ is free, so it follows that $S=(C \times F) / G$ is the desired surface.

\subsection{Surfaces of type $II$}
In this case $G=S_3$. \\
We denote by $r$ a $3-$cycle on $G$, and by $s$ a transposition.
Let $E$ be an elliptic curve, and let $p \in E$ be a point. We are looking for a smooth curve $C$ of genus $3$ such that there exists a $G-$cover $h: C \longrightarrow E$ branched at the point $p$. Hurwitz formula says that the ramification divisor has degree $4$, hence we have $h^{-1}(p)=q_1+q_2$ with $q_1,q_2 \in C$ distinct points, and the stabilizer group of $q_i$ is $\{ id, \; r, \; r^2 \} \; \cong  \mathbb{Z}_3$. The existence of such a double cover follows by the existence of the surjective homomorphism of groups:
\begin{displaymath}
\delta: \Delta(3) \longrightarrow G,
\end{displaymath}  
defined by $ \delta(a)=s$, \; $\delta(b)=r$. Indeed $\delta([a,b])=srsr^2=r$ which has order exactly $3$. Now observe that since $s$ acts without fixed points on $C$, we have an {\'e}tale double cover $C \longrightarrow C/\langle s \rangle$, where $C/\langle s \rangle$ is a smooth curve of genus $2$. Hence Theorem \ref{iperellbiell} shows that $C$ is both hyperelliptic and bielliptic. \\
Now we give an explicit equation for the curve $C$. By construction all the elements of order $2$ of $G$ act without fixed points, then the hyperelliptic involution  $\tau$ of $C$ does not belong to $G$. Hence $\widehat{G}=G$ and the action of $G$ on $C$ descends to the following action of $G$ on $\mathbb{P}^1$:
\begin{equation} \label{G2}
\begin{split}
r(x_0,x_1) &=(x_0, \; \xi x_1) \quad \textrm{where} \; \; \xi:=e^{\frac{2 \pi i}{3}} \\
s(x_0,x_1) & =(x_1, \; x_0).
\end{split}
\end{equation}
As in the previous example, we write $p(x)= \prod_{i=1}^8 (x-\alpha_i)$, and we notice that, since the zero locus $\Delta$ of $p$ is invariant under the action (\ref{G2}), the set $\{ \alpha_1, \ldots, \alpha_8 \}$ is of type $\{\alpha, \xi \alpha, \xi^2 \alpha, \frac{1}{\alpha}, \frac{1}{\xi \alpha}, \frac{1}{\xi^2 \alpha}, 0, \infty \}$. If we set $a:=\alpha ^3$, we get $h(x)=x(x^3-a)(x^3- \frac{1}{a})$, so the equation of $C$ is
\begin{equation} \label{C seconda}
y^2=x(x^3-a) \Big( x^3-\frac{1}{a} \Big) \quad a \neq 0, \; \pm 1,
\end{equation}   
and this proves that $C$ varies in one $1-$dimensional, irreducible family of curves. As in the previous case, there is only one lifting of the action of $G$ on $C$ with the desired stabilizers; it is the following:
\begin{displaymath}
\begin{split}
r(x,y) & =(\xi x, \xi^2 y) \\
s(x,y) & =(1/x,-y/x^4). 
\end{split}
\end{displaymath}
We can easily verify that $s$ acts without fixed points (because $a \neq \pm 1$), so the same holds for the other two transpositions of $G$. On the other hand, the fixed points of $r$ are $(0,0)$ and $(\infty,0)$, and their stabilizer group is $ \{ id, \; r, \; r^2 \} \cong \mathbb{Z}_3$. Since $s(0,0)=(\infty,0)$, these two points are in the same $G-$orbit, so $h:C \longrightarrow E$ is a $G-$cover branched at one point. \\
Now we can define the curve $F$ by the admissible epimorphism:
\begin{displaymath}
\mu: \Gamma(2,2,2,2,2,2) \longrightarrow G
\end{displaymath}
given by:
\begin{displaymath}
\begin{split}
x_1,x_2 & \longrightarrow s \\
x_3,x_4 & \longrightarrow rs \\
x_5,x_6 & \longrightarrow r^2s. \\
\end{split}
\end{displaymath}
In this way we constructed a $G-$cover $F$ of $\mathbb{P}^1$ branched
at $p_1, \ldots, p_6$ such that the $3-$cycles act without fixed
points; then Hurwitz formula gives $g(F)=4$. Now look at a generic triple cover of $\mathbb{P}^1$ branched in six points. A model of such a cover is the projection $\pi_p: C_3  \longrightarrow L$ of a smooth cubic curve $C_3 \subset \mathbb{P}^2$ onto a line $L$ from a general point $p \notin C_3$. It is classically known that from $p$ can be drawn six tangent lines to $C_3$, and since $p$ is general none of them is an inflexional tangent; whence $\pi_p$ has only simple ramification. Notice that  the $G-$cover $F$ is the galoisian closure of a cover of type $\pi_p$. Since the classic L$\ddot{\textrm{u}}$roth-Clebsch theorem states that the space of generic triple covers of $\mathbb{P}^1$ is connected (see [BaCa97, p.425]), it follows that $F$ belongs to a connected family (of dimension $3$). \\ 
Let $g_r=g(F/\langle r \rangle), \ g_1=g(F/\langle rs \rangle)$. Since
$S_3=D_3$, formula (\ref{accolaformula3a}) applies and we obtain
$4=g_r+2g_1$. On the other hand, since the $3-$cycles of $G$ act
without fixed points on $F$, then $F \longrightarrow F/\langle r
\rangle$ is an {\'e}tale $\mathbb{Z}_3-$cover, and the Hurwitz formula gives $6=3(2g_r-2)$, that is $g_r=2$. This in turn implies $g_1=1$, so $F \longrightarrow F/\langle sr \rangle$ is a double covering of an elliptic curve. Since $g(F)=4$, from Remark \ref{biell no 4} it follows that in this case $F$ is \emph{not hyperelliptic}. 
It is clear by our construction that the diagonal action of $G$ on $C \times F$ is free, and $S= (C \times F)/G$ is the desired surface. 

\subsection{Surfaces of type $III$}
In this case we have $G=D_4$.\\ 
We denote by $r$ a rotation and by $s$ a reflection. Let $E$ be an
elliptic curve, and let $p \in E$ be a point. We are looking for a
smooth curve $C$ of genus $3$ such that there exists a $G-$cover $f: C
\longrightarrow E$ branched at the point $p$. The Hurwitz formula
implies that the ramification divisor has degree $4$, hence we have
$f^{-1}(p)=q_1+q_2+q_3+q_4$, and the stabilizer group of each $q_i$ is
isomorphic to $\mathbb{Z}_2$. Now, notice that $r^2$ cannot act
without fixed points. Indeed, if this happens, then the same holds for
$r$ and $r^3$ (because $(r^3)^2=r^2$), and this in turn implies that
the $\mathbb{Z}_4-$cover $C \longrightarrow C/\langle r \rangle$ is
{\'e}tale. But if this is the case, then the Hurwitz formula gives
$2g(C/\langle r \rangle)-2=1$, which is impossible. Then $r^2$ has
fixed points, and since the stabilizer groups of points in the same
fibre are conjugated, this means that $\{ id, r^2 \}$ is the
stabilizer of $q_1, \ldots, q_4$ and that all the other elements of
$G$ act on $C$ without fixed points. In particular the reflection $s$
acts without fixed points, hence $C \longrightarrow C/\langle s
\rangle$ is an {\'e}tale double cover of a curve of genus $2$, so we again  
apply Theorem \ref{iperellbiell} in order to conclude that $C$ is both  hyperelliptic and bielliptic. \\      
Now we want to give an explicit equation for the curve $C$. Again $\tau \notin G$, then the action of $G$ on $C$ descends to the standard  action of $G$ on $\mathbb{P}^1$: 
\begin{equation} \label{G3}
\begin{split}
r(x_0,x_1)=(x_0, \; ix_1) \\
s(x_0, x_1) = (x_1, \; x_0).
\end{split}
\end{equation}
Let $p(x)=\prod_{i=1}^8 (x-\alpha_i)$; if we require that the zero locus $\Delta$ of $p$ is invariant under (\ref{G3}), we see that the set $\{ \alpha_1, \ldots ,\alpha_8 \}$ must be of the form $\{ \alpha, i\alpha, -\alpha, -i \alpha, \frac{1}{\alpha},  \frac{1}{i \alpha},  -\frac{1}{\alpha},  -\frac{1}{i \alpha} \}$. If we set $a:= \alpha ^4$, then the polynomial $p(x)$ will be of type $p(x)=(x^4-a)(x^4-\frac{1}{a})$, thus the equation of the curve $C$ is
\begin{equation} \label{C terza}
y^2=(x^4-a) \Big( x^4-\frac{1}{a} \Big) \quad  a \neq \pm 1,
\end{equation}   
and this in turn shows that $C$ belongs to an irreducible family of dimension $1$. Lift the action of G on $C$ as follows:
\begin{displaymath}
\begin{split}
r(x,y) & =(i x, -y) \\
s(x,y) & =(1/x,-y/x^4).
\end{split} 
\end{displaymath}
Notice that:
\begin{itemize}
\item[-] since $(0,0) \notin C$, $r$ acts without fixed points, so the same holds for its conjugate $r^3$;
\item[-] $a \neq 1$ implies that $s$ acts without fixed points, so the same holds for its conjugate $r^2s$;
\item[-] $a \neq -1$ implies that $rs$ acts without fixed points, so the same holds for its conjugate $r^3s$.
\end{itemize}
Moreover  $(0,1),(0,-1), (\infty, 1), (\infty, -1)$ are the only points on $C$ with non-trivial stabilizer, and this is $\{ id, r^2 \} \cong \mathbb{Z}_2$. All these four points are conjugated under the action of $G$, hence the $G-$cover $C \longrightarrow C/G$ is branched in just one point. As before, the other possible lifting does not give an action with the desired stabilizers. \\
Now we have to construct $F$. We know that the non-trivial stabilizer groups of points of $F$ are isomorphic to $\mathbb{Z}_2$. On the other hand, since we have shown that  $r^2$ has fixed points on the curve $C$ and we know that the diagonal action of $G$ must be free on the product $C \times F$, it follows that $r^2$ acts on $F$ without fixed points. Hence $\langle r \rangle=\{  id,r,r^2,r^3 \} \cong \mathbb{Z}_4$ acts without fixed points.
This means that the $G-$cover $f:F \longrightarrow \mathbb{P}^1$ is defined by an admissible epimorphism
\begin{displaymath} 
\mu: \Gamma(2,2,2,2,2,2) \longrightarrow G 
\end{displaymath}
such that the following condition is verified:
\begin{itemize}
\item[($**$)] the images of the generators $x_1, \ldots, x_6$ are elements of order $2$ different from $r^2$. 
\end{itemize}
Such a homomorphism exists; for instance we can consider the one given in the following way:
\begin{displaymath} 
\begin{split}
x_1, \ldots , x_4 & \longrightarrow s \\
x_5, \; x_6 & \longrightarrow rs.
\end{split} 
\end{displaymath} 
 Now recall that we have $5$ conjugacy classes in $G$:
\begin{displaymath} \{id \}, \; \{r^2\}, \;  \{s, r^2 s \}, \;  \{rs, r^3 s \}, \; \{r, r^3\}, \end{displaymath} and the center of the group is $Z:=\{ id, r^2 \} \cong
\mathbb{Z}_2$. $G$ is a central extension
\begin{displaymath} 
1 \longrightarrow Z \longrightarrow G \xrightarrow{\; \; \pi \;} (\mathbb{Z}_2)^2 \longrightarrow 1
\end{displaymath} 
such that two conjugated elements of $G$ are mapped by $\pi$ to the same element of $(\mathbb{Z}_2)^2$. Choose a basis $\{ e_1, \; e_2 \}$ of
$(\mathbb{Z}_2)^2$ such that the image of the conjugacy class $\{s, r^2s \}$ is $e_1$ and the image of the conjugacy class $\{ rs, r^3s \}$ is $e_2$. Let $\mu:\Gamma(2,2,2,2,2,2) \longrightarrow G$ be an admissible epimorphism
defining a $G-$cover $f:F \longrightarrow \mathbb{P}^1$, and let
$g_i=\mu(x_i)$ as usual. Moreover, let $n_j$ be the cardinality of the set $\{ i \; | \; \pi(g_i)=e_j \}$. Since the product $g_1g_2 \cdots g_6$ is equal to $1$, we have only two possibilities, that is either
$n_1=4, \; n_2=2$ or $n_1=2, \; n_2=4$. But it is clear that these two cases are $\textrm{Aut}(G)-$equivalent: indeed, the automorphism of $G$ given by $r \longrightarrow
r, \; s \longrightarrow rs$ interchange them. Then without loss of generality we can suppose that the former case occurs. If $x \in G, \; x \neq r^2$, denote by $x'$ the
element in the same conjugacy class of $x$, different from $x$. By a straightforward computation it is easy to see that:

\begin{enumerate}
\item[($1$)] a Hurwitz move on $G$ induces a permutation of the $\pi(g_i)$, since $(\mathbb{Z}_2)^2$ is abelian;
\item[($2$)] if $x, \; y$ belong to the same conjugacy class of $G$, then $xy=yx$;
\item[($3$)] if $x, \; y$ belong to different conjugacy classes, then
  a Hurwitz move sends $(x,y)$ to $(y', x)$, whereas the inverse of a
  Hurwitz move  sends $(x,y)$ to $(y,x')$;
\item[($4$)] if $x,y$ belong to different conjugacy classes, then
  twice applying either  a Hurwitz move or the inverse of a Hurwitz move we send $(x,y)$ to $(x',y')$.
\end{enumerate}
By $(1)$ we can suppose $g_1, \ldots, g_4 \in \{s, r^2s \}$ and $g_5,g_6 \in \{rs, r^3s \}$. Moreover, modulo an inner automorphism of $G$, we can suppose $g_1=s$.

\begin{claim} We can suppose $g_1=g_2=g_3=g_4 =s$.  \end{claim} Indeed if $g_1=g_2=g_3=s, \; g_4=r^2s$, applying ($4$) we transform the string $(s, \;  s, \; s, \; r^2s,
\; g_5, \; g_6)$ into $(s, \; s, \; s, \; s, \; g_5', \; g_6)$. Then suppose $g_1=g_2=s$, $ \; g_3=g_4=r^2s$. Again by $(4)$, applying twice a Hurwitz move we can
transform the string $(s, \; s, \; r^2s, \; r^2s, \;  g_5 , \; g_6)$ into the string $(s, \; s, \; r^2s, \; s , \; g_5' , \; g_6)$, then we apply $(3)$ and we send it into
the string $(s, \;  s , \; r^2s, \; g_5 , \; s , \; g_6)$, and again applying $(4)$ we can transform it into the string $(s , \; s , \; s , \; g_5', \; s, \; g_6)$;
finally we use $(3)$ and we obtain $(s , \; s , \; s , \; s , \; g_5 , \; g_6)$. This proves our claim. \\ \\ Now, since $g_1g_2 \cdots g_6=1$, it remains to prove that
with a sequence of moves of type either $(3)$ or $(4)$ it is possible to transform the string $(s , \; s, \;  s , \; s , \; rs, \; rs)$ into the string $(s, \; s, \; s, \;
s , \; r^3s , \; r^3s)$. But this is very easy: \begin{displaymath} \begin{split} ( s, \; s, \;  s, \; s, \; rs, \; rs) \xrightarrow{\; \; (3) \;} \;  (s , \; s, \;  s, \;
r^3s , \; s , \; rs ) \xrightarrow{\; \; (3) \;} \;  (s , \; s , \; s , \; r^3s , \; r^3s , \; s ) & \\ \xrightarrow{\; \; (3) \;} \; (s , \; s , \; s , \; r^3s , \; r^2s
, \; r^3s) \xrightarrow{\; \; (3) \;} \;  (s , \; s , \; s , \; s , \; r^3s , \; r^3s). \end{split} \end{displaymath} Hence all the admissible
epimorphisms $\mu$ which verify the condition ($**$) are Hurwitz equivalent, therefore Proposition \ref{componenti di Hurwitz} shows that $F$ belongs to a connected family of dimension $3$.\\
We know that the only elements of $G$ fixing points on $F$ are
$s,r^2s, rs, r^3s$. Notice that both $s$ and $r^2s$ fix $4$ points,
whereas both $rs$ and $r^3s$ fix $8$ points. If we set
$g_1=g(F/\langle rs \rangle), \  g_2=g(F /\langle r^2s \rangle)$, then
the Hurwitz formula yields $g_1=1, \; g_2=2$. Then in this case $F$ is a curve of genus $5$ which is a double cover of an elliptic curve, hence Remark \ref{biell no 4} implies that $F$ is \emph{not hyperelliptic}. By construction the action of $G$ on the product $C \times F$ is free, therefore $S=(C \times F) /G$ is a surface with the desired invariants.

\section{The bicanonical map} \label{bicanonica mappa}
In Section \ref{p_g=q=1} we showed that if $S$ is a surface with
$p_g=q=1$ and bicanonical map of degree $2$, then $S$ is a surface of
type $I, \; II$ or $III$. In Section \ref{esempi} we constructed these
surfaces, showing that they do exist. Now we will prove, conversely,
that \emph{any} surface isogenous to a product of type $I, \; II, \;
III$ has a bicanonical map of degree $2$, and that this map factors through the involution $\sigma$ induced on $S$ by $\tau \times id: C \times F \longrightarrow C \times F$, where $\tau$ as usual denotes the hyperelliptic involution of $C$. This means that $\sigma$ coincides with the bicanonical involution of $S$. Let us begin with some preliminary results.
   
\begin{lemma} \label{canoniconumerico}
If $S$ is a surface of type $I$, $II$, $III$, then $K_S$ is numerically equivalent to $C+4(K_SF)^{-1}F$.
\end{lemma}
\begin{proof}
This is an immediate consequence of the algebraic index theorem. Indeed for such a surface we have $(C+4(K_SF)^{-1}F-K_S)K_S=0$ and $(C+4(K_SF)^{-1}F-K_S)^2=0$.
\end{proof}

\begin{proposition} \label{canonicolineare}
Let $S$ be a surface as above. Then we have:
\begin{itemize}
\item $2K_S \cong 2C + 2\Phi_1+ 2\Phi_2$ \quad  if $S$ is of type $I$;
\item $2K_S \cong 2C + 4 \Phi_1$ \quad if $S$ is of type $II$;
\item $2K_S \cong 2C + 2 \Phi_1$ \quad  if $S$ is of type $III$.
\end{itemize}
\end{proposition}
\begin{proof}
This follows from Serrano's formula for the canonical class of a quasi-bundle (see [Se96, Theorem 4.1]).
\end{proof}

\begin{proposition}\label{gradobic=2}
Let $S$ be a surface of type $I$, $II$ or $III$. Then the bicanonical map of $S$ factors through the involution $\sigma$ induced on $S$ by $\tau \times id: C \times F \longrightarrow C \times F$.  
\end{proposition}
\begin{proof} 
Let $T:= S / \sigma$ and let $\psi: S \longrightarrow T$ be the projection onto the quotient. We denote by $q_{\alpha}, \; q_{\beta}$ the fibrations induced on $T$ by $\alpha: S \longrightarrow E$ and $\beta: S \longrightarrow \mathbb{P}^1$. Notice that the general fibre of $q_{\beta}$ is a smooth rational curve, because $C / \tau \cong \mathbb{P}^1$. Then $T$ contains a rational pencil whose general element is isomorphic to $\mathbb{P}^1$, and this implies that $T$ is a rational surface; in particular we have $p_g(T)=q(T)=0$. It follows that the inverse image in $S$ of a general fibre of $q_{\alpha}$ is not connected, otherwise $q(T)=1$. This means that the involution $\sigma$ induces a non-trivial involution on the elliptic curve $E=\textrm{Alb}(S)$, and this in turn implies that the divisorial fixed locus $R$ of $\sigma$ is composed of Albanese fibres; then $FR=0$. Moreover it is clear that $CR=8$, since a smooth curve of genus three contains $8$ Weierstrass points; therefore Lemma \ref{canoniconumerico} shows $K_S R=8$. Now we can compute the number $t$ of isolated fixed points of $\sigma$ by the holomorphic fixed point formula (see [DMP02, p.359]); we have
\begin{displaymath}
\sum_{i=0}^2 (-1)^i \textrm{Trace}(\sigma|H^i(S, \mathcal{O_S}))= \frac{t-K_S R}{4}.
\end{displaymath}
Since $T$ is rational and $p_g=q=1$, the left-hand side of this equation is equal to $1$, therefore we obtain $t=12$. So we have $t=K_S^2+4$, and Proposition \ref{fattorizzazione} allows us to conclude that $\phi$ is composed with $\sigma$.
\end{proof}
Proposition \ref{gradobic=2} says that there exists a commutative diagram
\begin{equation} \label{triang 6}
\begin{CD} 
\xymatrix{
S \ar[rrrr]^{\phi} \ar[rrd]_{\psi} & & & & \Sigma \\
& & T \ar[rru]_{\eta} & &
} 
\end{CD}
\end{equation}
By Proposition \ref{K_S^2 geq 5} we know that the degree of the bicanonical map is either $2$ or $4$, hence the degree of $\eta$ is either $1$ or $2$. 

\begin{proposition} \label{claim 3.1}
If the degree of $\phi$ is $4$, then the image in $\Sigma$ of the pencil $|C|$  is a base point free pencil of conics $|M|$.
\end{proposition}
\begin{proof}
By adjunction, the linear system $|2K_S|$ cuts out on $C$ a subseries of the bicanonical series $|2K_C|$. If we denote by $M$ the image of $C$ in $\Sigma$, by Proposition \ref{gradobic=2} we have two possibilities:
\begin{itemize}
\item[-] $\phi|_{C}$ has degree $2$ and $M$ is a quartic curve in $\mathbb{P}^4$;
\item[-] $\phi|_{C}$ has degree $4$ and $M$ is a conic.
\end{itemize}
If deg$\thinspace \phi=4$, then the second case occurs. Indeed, by Proposition \ref{canonicolineare} it follows that there exist bicanonical divisors of the form $C_1+C_2+\Phi$, where $C_1, \; C_2 \in |C|$ and $\Phi \neq 0$; they correspond  to hyperplane sections of $\Sigma$ of the form $M_1+M_2+\Psi$, where $M_1, \; M_2 \in |M|$ and $\Psi \neq 0$ (notice that the curve $\Phi$ cannot be contracted by $\phi$ because $K_S$ is ample). If $\textrm{deg}(\phi)=4$ then the degree of a hyperplane section of $\Sigma$ is $8$, and we have $2  \cdot$deg($M$)$=8-$deg($\Psi$)$<8$. This means $\textrm{deg}(M)<4$, that is $\textrm{deg}(M)=2$. Moreover $|M|$ is base point free, since $|C|$ is base point free and $\phi$ does not contract curves as $K_S$ is ample.
\end{proof}

\begin{remark}
Another way to prove Proposition \ref{claim 3.1} is the following: Proposition \ref{canonicolineare} tells us that $|2K_S|$ separates the curves of $|C|$, hence the degree of $\phi$ must be equal to the degree of the restriction of $\phi$ to $C$.
\end{remark}

\begin{proposition} \label{rgradobic=2}
If $S$ is of type $I, \; II, \; III$ then the degree of $\phi$ is $2$.
\end{proposition}
\begin{proof}
By contradiction, suppose that the degree of $\phi$ is $4$. Then $\Sigma$ is a linearly normal surface of degree $8$ in $\mathbb{P}^8$, and moreover it is rational because it is dominated by the rational surface $T$. Then it is well known (see [Na60, Theorem 8]) that $\Sigma$ is either the Veronese embedding in $\mathbb{P}^8$ of a quadric $Q \subset \mathbb{P}^3$, or the image of the blowup $\widehat{\mathbb{P}}$ of $ \mathbb{P}^2$ at a point $p$ via its anticanonical map $f:\widehat{\mathbb{P}} \hookrightarrow \mathbb{P}^8$. 
Now, observe that Proposition \ref{claim 3.1} implies that if deg $\phi=4$ then $\Sigma$ is a \emph{smooth} surface. Indeed, if $\Sigma$ is singular, it must be the Veronese embedding in $\mathbb{P}^8$ of a quadric cone $Q \subset \mathbb{P}^3$, and the linear system $|M|$ must be the system of conics through the node of $\Sigma$, which is impossible because $|M|$ is base point free. Now, let us consider separately the two possibilities for $\Sigma$. \\ \\
\emph{Case 1}.
$\Sigma$ is the Veronese embedding of a smooth quadric. Then $\Sigma$ contains two rulings of conics $|M|$, $|N|$, where $M^2=N^2=0,\; MN=1$. Proposition \ref{claim 3.1} allows us to suppose that $|M|$ is the image of the pencil $|C|$. Let $|D|$ be the inverse image in $S$ of the pencil $|N|$, and let $D \in |D|$ be a general curve. We have $K_SD=4,\;D^2=0$; therefore, since by Proposition \ref{caso standard} $S$ does not contain any genus $2$ pencils, it follows that $D$ is irreducible, that is, $|D|$ is a base-point free rational pencil of genus $3$ curves on $S$; moreover $MN=1$ gives $CD=4$. Applying Lemma \ref{canoniconumerico} we obtain
\begin{displaymath}
4=K_SD= \Big( C+\frac{4}{K_SF}F \Big) D=4+4\frac{FD}{K_SF},
\end{displaymath}
thus $FD=0$. But this is a contradiction, because $h^0(S, \;F)=1$,
whereas $h^0(S, \;D)=2$.\\ \\ 
\emph{Case 2}. 
$\Sigma$ is the embedding of $\widehat{\mathbb{P}}$ via its anticanonical system. Then $\Sigma$ contains exactly one pencil of conics $|M|$, which comes from the lines of $\mathbb{P}^2$ passing through the point that we have blown-up. Proposition \ref{claim 3.1} shows that $|M|$ is the image of the pencil $|C|$. We can exclude this case in two different ways.\\ \\ 
($1$) Let $E$ be the unique curve in $\Sigma$ such that $E^2=-1$, and
let $P:=\phi^*(E)$. Of course $P^2=-4$, and since $E \subset
\mathbb{P}^8$ is a line and $\textrm{deg} \; \phi=4$, we obtain
$K_SP=2$; it follows that $p_a(P)=0$. If $P$ is irreducible and
reduced, then $P$ is a smooth rational curve, and $P$ must then be
contained in some reducible Albanese fibre of $S$, a
contradiction. Therefore, let us suppose $P$ is reducible. Since $K_S$ is ample we can write $P=P_1+P_2$, with $P_i$ irreducible and $K_SP_i=1$; on the other hand, the index theorem gives $K_S^2P_i^2 \leq (K_SP_i)^2$, that is $P_i^2 \leq 0$; since $S$ does not contain rational curves, we obtain $P_1^2=P_2^2=-1$. Now $P^2=-4$ implies $P_1P_2=-1$, so  that $P$ is  non-reduced, and we have $P=2P'$, where $P'$ is reduced and irreducible, $K_SP'=1, \; (P')^2=-1$. Of course $CP'\geq 1, \; FP'\geq 1$, so applying Lemma \ref{canoniconumerico} we obtain:
\begin{displaymath}  
1=K_SP' = \Big( C+\frac{4}{K_SF}F  \Big) P'>1,
\end{displaymath}
a contradiction. \\  \\
($2$) Look at diagram (\ref{triang 6}). Since $K_S$ is ample, $\phi$
is a finite map; but $\psi$ is a finite map as well, hence $\eta$ must
be a finite map. As deg$\thinspace \phi=4$ and deg$ \thinspace
\psi=2$, it follows that deg$\thinspace \eta=2$, that is, $\eta: T \longrightarrow \Sigma$ is a double cover. Denote by $2aM+2bE$ the class of the branch divisor $B_{\eta} \subset \Sigma$ of $\eta$, where $a, \; b$ are integers, and $M, \; E$ are as above. Recall moreover that the canonical resolution $\hW$ of the singularities of $T$ is isomorphic to a Hirzebruch surface blown up in $12$ points, hence its invariants are $\chi(\mathcal{O}_{\widehat{W}})=1, \ K_{\widehat{W}}^2=-4$. Since $T$ has only ordinary double points as singularities, formulae (9) in [BPV, p.183] apply and we can write down the following equations:
\begin{equation} \label{invarianti di W}
\left \{ \begin{array}{ll}
1= 2 +\frac{1}{2}(-3M-2E)(aM+bE)+\frac{1}{2}(aM+bE)^2 \\
-4=16 + 4(-3M-2E)(aM+bE)+2(aM+bE)^2.
\end{array} \right. 
\end{equation}
A straightforward computation shows that the only pairs $(a, \; b)$ satisfying (\ref{invarianti di W}) are $(5/2, \; 3)$ and $(7/2, \;  1)$; this is absurd because $a$ and $b$ must be integer numbers. 
\end{proof}
Now, looking at Remarks \ref{vale rat inv} and \ref{vale rat inv 2} it is clear that we have actually proved the following:

\begin{theorem} \label{grado2fattor} 
Let $S$ be a minimal surface of general type with $p_g=q=1, \; K^2=8$ and let $\phi$ be the bicanonical map of $S$. Then the following are equivalent:  
\begin{itemize}
\item[$($a$)$] $\phi$ factors through a rational involution; 
\item[$($b$)$] $\phi$ has degree $2$;
\item[$($c$)$] $S$ is a surface of type $I, \; II$ or $III$.
\end{itemize}
\end{theorem}

\section{The plane models} \label{modellipiani}
The surfaces of type $I, \; II, \; III$ possess a rational involution $\sigma$ induced by the involution $\tau \times id: C \times F \longrightarrow C \times F$, where $\tau$ is the hyperelliptic involution of the curve $C$. We have shown in the last section that the involution $\sigma$ coincides with the bicanonical involution of $S$, hence it is possible to study the corresponding quotient map $S \longrightarrow S/ \sigma$ using the techniques introduced in Section \ref{degree 2}. A similar analysis in the case $p_g=q=0, \; K_S^2=8$ can be found in [Par03]. We remark that in the recent work [Bor03] Borrelli showed that if the bicanonical map  of a surface $S$ presenting the non-standard case factors through a double cover $\psi$ of a rational surface, then $S$ is regular unless $p_g=q=1$; he also described the plane model associated to $\psi$ in the last case, obtaining a weak version of our Theorem \ref{plane models}.\\ Let $\hW=\hS / \hat{\sigma}$; then $\hW$ is a Segre-Hirzebruch surface blown-up in $12$ points. Now let $W \cong \mathbb{F}_e$ be a good minimal model of $\hW$, and let $B$ be the branch locus of the generically double cover $S \longrightarrow W$. Then formula (\ref{espr B}) gives
\begin{displaymath}
B \cong 8C_0+(18+4e)L,
\end{displaymath}     
where $C_0$ is the section such that $C_0^2=-e$.
Proposition \ref{fibre in ram} implies that $B$ contains six fibres $L_1, \ldots , L_6$ corresponding to the six double curves of $|C|$; then we can write
\begin{displaymath} 
B=B^{\sharp}+L_1+ \cdots +L_6,
\end{displaymath}
where $B^{\sharp} \cong 8C_0+(12+4e)L$.
Over each curve $L_i$ there is a point $R_i$, and the curve $B$ contains a point $[5,5]$ at each point $R_i$, and no other singularities. Moreover, all these points $[5,5]$ are tangent to $L_i$ at $R_i$; therefore the only singularities of the curve $B^{\sharp}$ are six points of type $[4,4]$, one for each point $R_i$. This, together with the fact that these $[4,4]-$points have ``vertical'' tangent, implies that no section of $W$ is contained in $B^{\sharp}$. Thus we have $0 \leq C_0B^{\sharp}=12-4e$, that is $e \leq 3$ and at least $3+e$ of the points $R_i$ are not in $C_0$. If $e >1$, choose one of these points $R_i$ and perform an elementary transformation centered in it; in this way one can replace $e$ by $e-1$, and finally obtain $e=1$. Then we can contract the exceptional curve $C_0$ of $W=\mathbb{F}_1$ to a point, obtaining a birational map $f: S \longrightarrow \mathbb{P}^2$ which is generically a double cover. Following [Par03], we will call it a \emph{plane model} of the surface $S$. Let us denote again by $B$ the branch curve of this map, and let $P \in \mathbb{P}^2$ be the image of $C_0$.

\begin{theorem} \label{plane models}
Let $S$ be a minimal surface of general type with $p_g=q=1, \; K_S^2=8$ and bicanonical map of degree $2$. Let $f: S \longrightarrow \mathbb{P}^2$ be a plane model of $S$, and $B \subset \mathbb{P}^2$ be the branch curve. Therefore: \\ 
$1.$ $B=C_{16}+L_1+ \ldots +L_6$, where $C_{16}$ is a curve of degree $16$ and $L_1, \ldots, L_{6}$ are distinct lines passing through a point $P$. \\ 
$2.$ The singularities of $C_{16}$ are:
\begin{itemize}
\item a singular point of multiplicity $8$ at $P$;
\item six points $\xi_i:=(x_i,y_i)$ of type $[4,4]$, such that $\xi_i$ is tangent to $L_i$ at a point $R_i$. 
\end{itemize}      
$3.$ The curve $C_{16}$ looks as follows:
\begin{itemize}
\item if $S$ is of type $I$, then $C_{16}=C_8^1+C_8^2$, where $C_8^1, \; C_8^2$ are irreducible curves of degree $8$ having a $4-$ple point at $P$ and a $[2,2]-$singularity (ordinary tacnode) tangent to $L_i$ at the point $R_i$, for any $i$; 
\item if $S$ is of type $II$, then $C_{16}=C_4+C_{12}$, where $C_k$ is an irreducible curve of degree $k$. Moreover $C_4$ has a double point at $P$ and it is tangent to $L_i$ at $R_i$, whereas $C_{12}$ has a $6-$ple point at $P$ and it has a $[3,3]-$point tangent to $L_i$ at $R_i$, for any $i$;
\item if $S$ is of type $III$, then $C_{16}$ is irreducible.
\end{itemize}
$4.$ There is exactly one conic $\Lambda$ containing the points $R_1, \ldots, R_6$. \\
Vice versa, the minimal model $S$ of the canonical resolution of a double cover  $X \longrightarrow \mathbb{P}^2$ branched along a curve $B$ as above is a surface with $p_g=q=1, \; K_S^2=8$ and bicanonical map of degree $2$.   
\end{theorem}
\begin{proof}
The proof of ($1$) and ($2$) is exactly as in [Par03] and it is left to the reader. The proof of ($3$) comes from ($1$), ($2$) and the description of the fixed locus of the bicanonical involution given in Remark \ref{fissodivisoriale}.
Now we give the proof of ($4$). Let $S^*$ be the canonical resolution of the double cover of $\mathbb{P}^2$ branched along the curve $B$; the invariants $K_{S^*}^2, \; \chi(\mathcal{O}_{S^*})$ of $S^*$ can be recovered from the degree of $B$, its singularities and their mutual position, using standard formulae: see for instance [BPV84, p.183]. In particular it follows that $p_g(S^*)=p_g(S)$ is equal to the dimension of the subspace of $H^0(\mathbb{P}^2, \; \mathcal{O}_{\mathbb{P}^2}(8))$ consisting of those sections vanishing of order at least $6$ in $P$, and such that their zero loci are curves with a point of multiplicity at least $2$ in $R_i$, and that are tangent to $L_i$ at $R_i$ for any $i$. If $C_8$ is such a curve, it follows that the intersection multiplicity of $C_8$ and $L_i$ at $R_i$ is at least $9$, hence $C_8$ splits in a curve of type $L_1+ \cdots +L_6+\Lambda$, where $\Lambda$ is a conic containing $R_1, \ldots, R_6$. Since $p_g(S)=1$, it follows that there is exactly one such a conic $\Lambda$. \\
Suppose now that we have a double cover $X \longrightarrow
\mathbb{P}^2$ branched along a curve $B$ as above, and let $S^*$ be
its canonical resolution. Straightforward computations show that the
minimal model of $S^*$ is a surface $S$ with $p_g=q=1, \;
K_S^2=8$. The pullback  of the pencil of lines through $P$ gives a
base point free genus $3$ pencil $|C|$ on $S$, which contains six
double fibres since $B$ contains the six lines $L_1, \ldots,
L_6$. Therefore the Zeuthen-Segre formula shows that the six double fibres are twice a smooth curve of genus $2$, whereas any other fibre of $|C|$ is smooth. Now a base-change argument as in the proof of Proposition \ref{isotrivial pencil} shows that $S$ is isogenous to a product, hence $S$ is a surface of type $I, \; II$ or $III$. Therefore we can apply Proposition \ref{gradobic=2} and we conclude that the degree of the bicanonical map of $S$ is $2$.
\end{proof}

\begin{remark}
The last part of Theorem \ref{plane models} suggests that one might
give another construction of surfaces of type $I, \; II, \; III$ if he
would be able to show directly the existence of the curves $B$. The
problem is that we ask for so many singularities in $B$ that these
cannot impose independent conditions on the linear parameters: in
other words, the linear system in $\mathbb{P}^2$ containing $B$ with
the prescribed singularities is very superabundant. In fact, its
expected dimension is negative. On the other hand, since we have shown
that surfaces of type $I, \; II, \; III$ do exist, it follows that the
curves $B$ exist as well. But since their effective construction seems
very difficult, we will not attempt to do it. This is an interesting
example of the fight ``\emph{Campedelli versus Godeaux}'', see [Reid91].
\end{remark}

\section{The moduli space} \label{mmoduli}

Let $ \mathcal{M}$ be the moduli space of minimal surfaces of general type with $p_g=q=1,  \; K_S^2=8$, and let $\mathcal{S}$ be the subset of $\mathcal{M}$ corresponding to the surfaces having the bicanonical map of degree $2$.  Let moreover  $\mathcal{S}_I, \; \mathcal{S}_{II}, \mathcal{S}_{III}$ be the subsets of $\mathcal{M}$ corresponding to surfaces of type $I$, $II$, $III$, respectively.

\begin{theorem} \label{centrale moduli}
The following hold:
\begin{enumerate}
\item[$(i \;)$] $\mathcal{S}$ is a disjoint union: $\mathcal{S}=\mathcal{S}_I \sqcup \mathcal{S}_{II} \sqcup \mathcal{S}_{III}$;
\item[$(ii \;)$] $\mathcal{S}_I, \, \mathcal{S}_{II}, \, \mathcal{S}_{III}$ are irreducible components of $\mathcal{M}$ of the following dimensions:
$$dim \; \mathcal{S}_I=5, \quad dim \; \mathcal{S}_{II}=4, \quad dim \;  \mathcal{S}_{III}=4;$$
\item[$(iii \;)$] $\mathcal{S}_I, \, \mathcal{S}_{II}, \, \mathcal{S}_{III}$ are normal varieties.
\end{enumerate} 
\end{theorem}
The proof of Theorem \ref{centrale moduli} will be a consequence of the following lemmas.

\begin{lemma} \label{sottospazi invarianti}
Let $X$ be a smooth variety, and let $G$ be a group acting faithfully on $X$ in such a way that the quotient $Y:=X/G$ is smooth. Let $\pi: X \longrightarrow Y$ be the quotient map and $D$ the reduced branch divisor. Then: 
\begin{itemize}
\item[] $(\pi_* \omega_X)^G=\omega_Y$;
\item[] $(\pi_* \omega_X^2)^G=\omega_Y^2(D)$.
\end{itemize}
\end{lemma}
\begin{proof}
See [Par91, Proposition 4.1] and [Par03, Lemma 6.4].
\end{proof}

\begin{lemma} \label{insiemi aperti}
$\mathcal{S}_I, \; \mathcal{S}_{II}, \mathcal{S}_{III}$ are open and normal subsets of $\mathcal{M}$.
\end{lemma}

\begin{proof}
We use a standard argument in deformation theory. For any variety $X$
let $\textrm{Def}(X)$ be the deformation functor of $X$; if there is on $X$ the action of a group $G$, let $\textrm{Def}(X, \; G)$ be the deformations of $X$ that preserve the $G-$action. Now, let $S$ be a surface of type $\mathcal{S}_I$,  [resp. $\mathcal{S}_{II}$, $\mathcal{S}_{III}$], and write $S=(C \times F) /G$. Roughly speaking, if we  deform $C$ and $F$ in such a way that $C$ remains hyperelliptic and the actions of $G$ are preserved, we obtain a deformation of $S$ having the same quotient structure as $S$; then Lemma \ref{insiemi aperti} will be proved if we show that the general deformation of $S$ is of this type. In order to do this we show some intermediate results.

\begin{claim} \label{lemma aperti}
Let $S$ be a surface of type $I$ [resp. $II$, $III$], $S=(C \times
F)/G$. Let $G_0$ be the subgroup of Aut($C$) generated by $G$ and
$\tau$, let $D$ be the branch locus of the $G-$covering $\;h: C
\longrightarrow E$ and let $D_0$ be the branch locus of the $G_0-$covering $\;h_0: C \longrightarrow \mathbb{P}^1$. Then $\mathcal{S}_I$ [resp. $\mathcal{S}_{II}$, $\mathcal{S}_{III}$] is an open and normal subset of $\mathcal{M}$ if $deg(D_0)=deg(D)+3$.
\end{claim}
Indeed, let $\eta$ be the natural map of functors $\eta: \textrm{Def}(C, \; G_0) \times \textrm{Def}(F, \; G) \longrightarrow \textrm{Def}(S)$. Since $C$ and $F$ are curves, their deformation spaces are unobstructed. On the other hand, the tangent space in $[S]$ to $\textrm{Def}(S)$ is given by
\begin{displaymath}
T_{[S]}\textrm{Def}(S)=H^1(S,\; T_S)=H^1(C, \; T_C)^G \oplus H^1(F, \; T_F )^G,\end{displaymath}
by Kunneth's formula and because $G$ acts separately on $C$ and $F$. So the differential of $\eta$ is simply the inclusion
\begin{displaymath}
\eta_*: H^1(C, \; T_C)^{G_0} \oplus H^1(F, \; T_F)^G \hookrightarrow H^1(C, \; T_C)^G \oplus H^1(F, \; T_F)^G.
\end{displaymath}
Notice that if $\eta_{*}$ is an isomorphism, then the Kuranishi family of $S$ has the same dimension as $\textrm{Def}(C, \; G_0) \oplus \textrm{Def}(F,\; G)$, and this in turn shows that $\mathcal{S}_I$ [resp. $\mathcal{S}_{II}$, $\mathcal{S}_{III}$] is an open subset of $\mathcal{M}$. Thus we are done if we prove that $H^1(C, \; T_C)^{G_0}$ is isomorphic to $H^1(C, \; T_C)^G$. By Serre duality, this is equivalent to proving that $H^0(C, \; 2 \omega_C)^{G_0}$ is isomorphic to $H^0(C, \; 2 \omega_C)^G$ or, by Lemma \ref{sottospazi invarianti}, that  $H^0(\mathbb{P}^1, \; 2K_{\mathbb{P}^1}+ D_0 )$ is isomorphic to $H^0(E, \; D)$. But this happens exactly when $\textrm{deg}(D_0)=\textrm{deg}(D)+3$. Notice that our argument shows as well that the Kuranishi family of $S$ is smooth; so the moduli space $\mathcal{S}$ is normal, and in particular it is the disjoint union of its irreducible components. 
 
\begin{claim} \label{splitting}
In all our examples, $G_0 \cong G \times \mathbb{Z}_2$.
\end{claim}
We have already remarked that in any case the hyperelliptic involution $\tau$ of $C$ does not belong to $G$. Thus, by Proposition \ref{carattere banale}, it is sufficient to prove that the character $\lambda$  defined in Subsection \ref{hyp aut} is trivial. We show that this is the case only when $S$ is of type $I$; in the other two cases, the proof is exactly the same. Suppose therefore that $S$ is a surface of type $I$, and let $S=(C \times F)/G$, where $G=\mathbb{Z}_2 \times  \mathbb{Z}_2$. Results of  Section \ref{surfaces of type I} imply that $C$ is isomorphic to the curve given in $\mathbb{P}(1,1,4)$ by the equation $y^2=p(x_0,x_1)$, where
\begin{displaymath}
p(x_0,x_1)=(x_1^2-ax_0^2)(x_1^2-bx_0^2) \Big( x_1^2-\frac{1}{a}x_0^2 \Big) \Big( x_1^2- \frac{1}{b}x_0^2 \Big).  
\end{displaymath}
The subgroup $H \subset SL(2,\mathbb{C})$ is generated by the two elements
\begin{displaymath}
h_1=\left( \begin{array}{cc}
-i & 0\\
0  & i \\
\end{array} \right),  \quad
h_2=\left( \begin{array}{cc}
0 & i\\
i & 0 \\
\end{array} \right);
\end{displaymath}
since both $h_1$ and $h_2$ act as the identity on $p(x_0, x_1)$, it follows that the character $\lambda$ is trivial. \\ \\
Now we can conclude the proof of Lemma \ref{insiemi aperti} by
separately analyzing the cases occurring for $S$. \\ \\
\emph{S is of type I.}\\
In this case the $G-$cover $h:C \longrightarrow E$ is branched at two
points, hence $\textrm{deg}(D)=2$.  Since by Claim \ref{splitting} we
have $G_0= G \times \mathbb{Z}_2=(\mathbb{Z}_2)^3$, each non-trivial
stabilizer in the $G_0-$covering $h_0:C \longrightarrow \mathbb{P}^1$
is isomorphic to $\mathbb{Z}_2$ and the Hurwitz formula yields $\textrm{deg}(D_0)=5$. Then we can apply Claim \ref{lemma aperti} and we are done.\\ \\   
\emph{S is either of type II or of type III.}\\ 
In this case the $G-$cover $h:C \longrightarrow E$ is branched at one point $p$, hence $\textrm{deg}(D)=1$. Moreover $G_0=G \times \mathbb{Z}_2$, and we have a commutative diagram
\begin{equation} \label{moduli diag II}
\begin{CD} 
\xymatrix{
C \ar[rrrr]^{h_0} \ar[rrd]_{h} & & & & \mathbb{P}^1 \\
& & E \ar[rru]_{g} & &
} 
\end{CD}
\end{equation}
where $g$ is a double covering. Let $q=g(p)$. We claim that $g$ is
branched at $q$. Suppose by contradiction that it is not. Then we have
$g^*(q)=p+p'$, where $p' \neq p$. Since $p$ is the only branch point
of $h$, it follows that $h^*(p+p')$ contains both points where $h$ is
ramified and points where $h$ is not ramified. This in turn implies
that the same holds for $h_0^*(q)$, a contradiction because $h_0$ is a
Galois covering. Then $g$ is branched at $q$, so the branch locus of
$h_0$ coincides with the branch locus of $g$, and this means
$\textrm{deg}(D_0)=4$. Then we can again use Claim \ref{lemma aperti}
in order to reach a conclusion. \\ \\
This ends the proof of Lemma \ref{insiemi aperti}.      
\end{proof}

\begin{remark} \label{deform hyp}
Notice that $H^0(C, \; T_C)^G=H^0(C, \; T_C)^{G_0}$ means that any deformation of the hyperelliptic curve $C$ which preserves the $G-$action is again hyperelliptic.  
\end{remark}

\begin{lemma} \label{chiusi}
$\mathcal{S}$ is a  closed subset of $\mathcal{M}$.
\end{lemma}
\begin{proof}
Let $\Delta \subset \mathbb{C}$ be the unitary disk and let $p \colon\mathcal{X} \longrightarrow \Delta$ be a smooth family of surfaces such that $X_t:=p^{-1}(t) \in \mathcal{S}$ for any $t \in \Delta-\{0\}$. We have to prove that $X_0$ belongs to $\mathcal{S}$ as well. Let $\sigma_t \colon X_t \longrightarrow X_t$ be the bicanonical involution of $X_t$; since $K_{X_t}$ is ample for any $t \in \Delta-\{0\}$, there is a birational map $\sigma \colon \mathcal{X} \dashrightarrow \mathcal{X}$ that restricts to $\sigma_t$ on $X_t$ for $t \neq 0$. By [FP97, Corollary 4.5] $\sigma$ is actually biregular and we denote by $\sigma_0$ the restriction of $\sigma$ to $X_0$. Since $\sigma_t$ has $12$ fixed points for any $t \neq 0$, we can use Cartan's lemma exactly as in [Par03, Lemma 6.3] in order to conclude that the involution $\sigma_0$ has $12$ fixed points, too. Moreover, since $X_t /\sigma_t$ is a rational surface for $t \neq 0$, it follows that $X_0 / \sigma_0$ is rational as well. Therefore Proposition \ref{fattorizzazione} tells us that the bicanonical map of $X_0$ factors through the rational involution $\sigma_0$, and this in turn implies by Theorem \ref{grado2fattor} that $X_0$ belongs to $\mathcal{S}$.
\end{proof}
\begin{lemma} \label{D irriducibili}
$\mathcal{S}_{I}, \; \mathcal{S}_{II}, \; \mathcal{S}_{III}$ are irreducible components of $\mathcal{M}$, whose dimensions are
\begin{displaymath}
\dim \mathcal{S}_{I}=5, \quad \dim \mathcal{S}_{II}=4, \quad \dim \mathcal{S}_{III}=4. 
\end{displaymath}
\end{lemma}

\begin{proof}
In Section \ref{esempi} we showed that in any case $C$ belongs to an
irreducible family of curves and $F$ belongs to a connected family. On
the other hand Claim \ref{lemma aperti} and Lemma \ref{chiusi} show that any deformation of $S$ comes from a deformation of $C$ and $F$ which preserves the $G-$action, therefore $\mathcal{S}_{I}, \; \mathcal{S}_{II}, \; \mathcal{S}_{III}$ are connected subsets of $\mathcal{M}$. Besides, since they are open and closed, and since $\mathcal{M}$ is normal, they are irreducible components of the moduli space. \\
We proved moreover that the curve $F$ varies in any case in a family of dimension $3$, whereas $C$ belongs to a family of dimension $2$ if $S$ is of type $I$ and to a family of dimension 1 otherwise. Therefore we have $\dim \mathcal{S}_{I}=5, \; \dim \mathcal{S}_{II}=4, \; \dim \mathcal{S}_{III}=4$. 
\end{proof}  
This completes the proof of Theorem \ref{centrale moduli}.

\section{Open problems} \label{open problems}

\noindent $\mathbf{Problem \; 1.}$ Construct (if they exist) surfaces of general type with $p_g=q=1, \; K^2=8$ which are not isogenous to a product. \\  \\
$\mathbf{Problem \; 2.}$ Classify all minimal surfaces of general type with $p_g=q=1$ which are isogenous to a product. In this case, we have necessarily $K_S^2=8$ (see [Se90]). We hope to deal with this topic in a forthcoming paper. The corresponding problem for surfaces with $p_g=q=0$ was recently considered by Bauer and Catanese; see [BaCa03]. \\ \\
$\mathbf{Problem \; 3.}$ Classify (if they exist) the minimal surfaces of general type with $p_g=q=1, \; K^2=8$ and bicanonical map of degree $4$. \\ \\

\end{document}